\documentclass[12pt]{article}
\usepackage[utf8]{inputenc}
\usepackage[T2A]{fontenc}

\usepackage[english]{babel}
\usepackage{fullpage}

\usepackage{amsmath}
\usepackage{mathrsfs} 
\usepackage{amsfonts}
\usepackage{amsthm}
\usepackage{bbm}
\usepackage{xcolor}
\usepackage{hyperref}
\usepackage{amssymb}
\usepackage{bm}
\usepackage{cite}
\usepackage{comment}
\usepackage{graphicx}
\usepackage{epstopdf}
\usepackage{algorithm}
\usepackage[noend]{algorithmic}

\usepackage{kamzolov}

\theoremstyle{definition}

\theoremstyle{plain}

\theoremstyle{remark}

\DeclareMathOperator*{\argminset}{Arg\,\min}

\newcommand{\y}{\mathtt{y}}
\newcommand{\p}{\mathtt{p}}
\newcommand{\z}{\mathtt{z}}
\newcommand{\w}{\mathtt{w}}

\def\e{\varepsilon}

\newcommand{\norm}[1]{\|#1\|}

\newcommand\uprule{\rule{0mm}{1.9ex}}

\newtheorem{Lm}{Lemma}%[section]

\newtheorem{Th}{Theorem}%[section]

{
\theoremstyle{remark}
\newtheorem{Rem}{Remark}%[section]

}
{
\theoremstyle{definition}
}

\begin{document}

\title{Composite optimization for the resource allocation problem}
 
\author{
Anastasiya Ivanova\textsuperscript{a}$^\ast$,
Pavel Dvurechensky\textsuperscript{b}$^{\dag}$,\\ Alexander Gasnikov\textsuperscript{a,c}$^{\ddag}$ and Dmitry Kamzolov \textsuperscript{a}$^{\S}$.\\ \\
\textsuperscript{a}Moscow Institute of Physics and Technology, Moscow, Russia;\\
\textsuperscript{b} Weierstrass Institute for Applied Analysis and Stochastics, Berlin, Germany;\\ 
\textsuperscript{c}Institute for Information Transmission Problems RAS, Moscow, Russia;
}

\maketitle

\begin{abstract}
In this paper we consider resource allocation problem stated as a convex minimization problem with linear constraints. To solve this problem, we use gradient and accelerated gradient descent applied to the dual problem and prove the convergence rate both for the primal iterates and the dual iterates. We obtain faster convergence rates than the ones known in the literature.  We also provide economic interpretation for these two methods. This means that iterations of the algorithms naturally correspond to the process of price and production adjustment in order to obtain the desired production volume in the economy. Overall, we show how these actions of the economic agents lead the whole system to the equilibrium.
\end{abstract}

\section{Introduction}

In this paper we consider a resource allocation problem in an economy consisting of distributed set of producers which are managed by a centralized price adjustment mechanism. Our approach is based on the state-of-the-art convex optimization methods, i.e. we consider the resource allocation as a convex optimization problem, solve it by first-order methods, provide convergence analysis, and give an economic interpretation of the steps of these methods. 

The problem of optimal resource allocation is to maximize producers’ aggregated profits by sharing available resources. Popularized and advocated mainly in the monograph \cite{Arrow}, the mechanisms of decentralized resource allocation gained a lot of attention in economics and operations research since then, see e.g. \cite{Friedman,Campbell,Kakhbod} and references therein. Each producer seeks to minimize its own costs and, in total, all the producers need to produce a certain amount of products. This problem can be cast as an optimization problem with the objective corresponding to the aggregated cost function of all producers and constraints corresponding to the condition for the necessary volume of production.
We assume that constraints are linear and separable. In this optimization problem, primal variables are production bundles and dual variables represent prices of resources.

Solution of this optimization problem depends on the way, how agents in this economy can exchange information with each other. One of the mechanisms is called Walrasian in the literature~\cite{Aubin}. This mechanism is an iterative process. At each iteration, some agents formulate a supply for others. Further, the second agents, based on the supply, form the optimal demand. Then, the first agents again form the optimal supply based on demand and so on until the system comes to Walrasian equilibrium. In the economic literature, this iterative process is called the Walrasian tâtonnement. For simplicity, let us consider an economy with only one good. We assume that there are some small producers and one big consumer (Center). Each producer has its own cost function and this functions are unknown to the Center. The Center sets a price for the good and producers provide the produced quantity for this price. The goal of the Center is to iteratively find such a price that all the factories produce the certain amount of good in total.  This will be a Walrasian equilibrium in the system. 
One example of such situation can be a centrally governed economic like in USSR. The price set by the Center plays the role of control parameter of the system and the goal is to provide production determined by plan. Another example could be farmers as producers and some large retail company as the Center. Retail company knows the demand for the good and sets the price so that the produced amount is equal to the demand. 
In the literature, an iterative numerical algorithm, which corresponds to Walrasian mechanism is based either on the dichotomy method (for one-product economy)~\cite{Ivanova} or the ellipsoid method (for many-product economy)~\cite{Friedman}. 
These algorithms are effective in the case when the amount of producers are small.

In this paper, following~\cite{NestShik}, we consider a different price adjustment mechanism.  In particular, we consider the resource allocation problem without centralized price control, each producer setting up its own price for selling products to the Center. The Center knows the amount which needs to be produced by all the producers in total, selects the most advantageous offers (i.e. selects the offers with the best price) and tries to purchase the product in the required volume. Producers adjust the volume of product and the prices, based on the volume bought from them by the Center and the demand from the Center for this particular factory. In their paper \cite{NestShik}, the authors use dual subgradient method with averaging as a numerical algorithm for this problem. The main advantage of their algorithm is that they provide convergence rate for the whole primal-dual sequence unlike the optimization literature which gives convergence rate for the running average. Their algorithm has optimal for convex nonsmooth optimization convergence rate $O\left(\dfrac{1}{\sqrt{t}}\right)$, $t$ being the iteration counter. 
In contrast to~\cite{NestShik}, in this paper, under an additional assumption of strongly convexity of the primal objective, we consider the dual optimization problem as a composite minimization problem, meaning that the objective in the dual problem is a sum of two functions a smooth and a simple non-smooth.
We use gradient descent 
to obtain convergence rate $O\left(\dfrac{1}{t}\right)$ and accelerated gradient descent to obtain convergence rate $O\left(\dfrac{1}{t^2}\right)$. 

The paper is organized as follows. In Section~2 we consider the primal problem and describe its economical interpretation. In Section~3 we describe the method of subgradient projection for the resource allocation problem and give the interpretation for the step. In Section~4 we use composite gradient method for resource allocation problem and obtain estimation for the convergence rate. In Section~5 we consider accelerated composite gradient descent. And in Section~6 we show some experiments, that verify our theory.

\section{Problem statement}
\label{Prob_stat}

In this section, we provide the statement of the resource allocation problem. For simplicity, we start with a one-product economy. The case of many products will be considered in Appendix. We assume that there is a Center and $n$ producers which produce one product. Each producer has its own cost function $f_k(x_k), \, k=1,\, \ldots \,,n$ representing the total cost of production of a volume $x_k \in \R$ -- the volume of product produced by the producer $k$ in one year. Since the producers are independent, the cost functions of the producers are unknown to the Center, and each producer knows only its own cost function. Each producer is also entitled to set its own price for product, but the price does not affect the quality of the product, i.e. all producers produce the same product, only at different prices. The Center buys product from the producers and chooses its strategy in such a way that the total production volume  per year by all producers is not less than $C$. To do so, the Center needs to find $y_k$ - the volume of product which is purchased from the producer $k$.  Then, each producer produces at least the volume $y_k$ of the product, and the goal of the system is to minimize the total cost of production. Thus, we consider the following resource allocation problem 
\begin{equation*}
\label{Main_problem}
(P) \quad \quad  \min\limits_{\substack{\sum \limits_{k=1}^n y_k \geqslant  C, \;x_k \geqslant  y_k;\; \\ y_k \geqslant  0,\;  x_k \geqslant  0,\; k = 1, \, \ldots, \, n,}}f(\mathtt{x})=\sum \limits_{k=1}^n f_k(x_k),
\end{equation*}
where cost functions~$f_k(x_k)\;,\ k = 1, \, \ldots, \, n$ are increasing and~$\mu$-strongly convex, i.e. \\ $f_k''(x) \geqslant \mu,  \,\forall x \geqslant 0, \,  \, k = 1, \, \ldots, \, n.$ 
\begin{Rem}
We point that, from the optimization point of view the variables $y_k$, $k=1,...,n$ are not necessary. We introduce them in order to the constructed dual problem and primal-dual method for solving the primal-dual pair of problems has an economic interpretation. The reformulation (P) is one of our contributions.
\end{Rem}

\begin{Rem}
Assumption of strong convexity of functions $f_k(x_k), \, k=1,\, \ldots \,,n$ holds, for example, when these functions are twice continuously differentiable and have positive second derivative. Economically this means that the production cost grows faster than the volume of the production. In other words, the production cost of a new unit of volume grows as the volume of production grows. For example, this happens for Agriculture. If the producer grows wheat, then the more he wants to produce from one hectare, the more he should invest in fertilizers, chemicals from pests, or even genetic technology. 
For a factory the producer has to invest more and more in more advanced facilities such as robots, production machines, etc.
\end{Rem}

Introducing dual variables $p_k$, $k=1,...,n$ and using the duality theory, we obtain 
\begin{equation*}
\begin{split}
\min \limits_{\substack{\sum \limits_{k=1}^n y_k \geqslant  C, \; x_k \geqslant  y_k,\; y_k \geqslant  0;\;\\ x_k \geqslant  0,\;k=1, \, \ldots, \, n}} f(\mathtt{x}) &= \min \limits_{\substack{\sum \limits_{k=1}^n y_k \geqslant  C\; y_k \geqslant  0;\;\\ x_k \geqslant  0,\;k=1, \, \ldots, \, n}} \Bigl\{f(\mathtt{x})+\sum \limits_{k=1}^n \max \limits_{p_k \geqslant  0}p_k(y_k-x_k)\Bigr\}\\
&=-\min\limits_{p_1, \, \ldots, \, p_n\geqslant  0} \Bigl\{\sum \limits_{k=1}^n \max \limits_{x_k \geqslant  0}(p_{k}x_k-f_k(x_k))-\min \limits_{\sum \limits_{k=1}^n y_k \geqslant  C; \;  y_k \geqslant  0}\sum \limits_{k=1}^{n}p_ky_k\Bigr\}\\
&[y^*_{k^*} = C, \;\text{where} \; k^* = \arg\min_kp_k \; \text{and} \; y_k^*=0 \; \text{for}\; k\ne k^*]\\
&=-\min\limits_{p_1, \, \ldots, \, p_n\geqslant  0}\Bigl\{\sum \limits_{k=1}^n \max \limits_{x_k \geqslant  0}(p_kx_k-f_k(x_k))-C\min \limits_{k=1, \, \ldots, \, n}p_k\Bigr\}\\
&=-\min\limits_{p_1, \, \ldots, \, p_n\geqslant  0} \Bigl\{\sum \limits_{k=1}^n \Big\{ p_kx_k(p_k)-f_k(x_k(p_k))\Big\}-C\min \limits_{k=1, \, \ldots, \, n}p_k\Bigr\},
\end{split}
\end{equation*}
where 
\begin{equation}
\label{x_k(p_k)}
x_k(p_k)=\argmax \limits_{x_k \geqslant  0} \Big \{ p_k x_k-f_k(x_k) \Big \}, \quad  k =1, \, 2, \, \ldots, \, n.
\end{equation}
Then the dual problem (up to a sign) has the following form
\begin{equation*}
\label{Dual_problem}
(D) \quad \quad \varphi(p_1, \, \ldots, \, p_n)=\sum \limits_{k=1}^n \Big\{ p_kx_k(p_k)-f_k(x_k(p_k))\Big\}-C\min \limits_{k=1, \, \ldots, \, n}p_k \rightarrow \min \limits_{p_1, \, \ldots, \, p_n\geqslant  0.}
\end{equation*}

Note that, the Slater's constraint qualification condition holds for the primal problem $(P)$. Thus, the strong duality holds and both the primal problem $(P)$ and the dual problem $(D)$ have solutions. Throughout the paper, we solve the dual problem by different first-order methods, interpret the steps of these methods and show, how the primal variables $x_k,y_k$, $k=1,...,n$ can be reconstructed. We point that primal-dual gradient methods \cite{dvurechensky2018decentralize,guminov2019accelerated,dvurechensky2016primal-dual} are not applicable here since the dual problem (D) is a composite optimization problem.

\section{Subgradient descent}
\label{Subgradient descent}
For the sake of completeness, in this section, we consider dual problem $(D)$ as a non-smooth optimization problem and apply subgradient method to solve it with the rate $O(1/\sqrt{t})$. We also provide an economic interpretation of the numerical procedure based on subgradient method. The material of this section is not new and mostly follows~\cite{Ivanova},
but we include it to be able to compare, in the next sections, its convergence rate and interpretation with faster approaches based on composite gradient descent.

To  solve the problem $(D)$, we use the projected subgradient method with the step given by 
\begin{equation}
\label{5_eq_subgrad_der_1}
\mathtt{p}^{t+1}= \left [ \mathtt{p}^{t}-h g(\mathtt{p}^{t}) \right ]_{+},
\end{equation} 
where $h$ is the stepsize, which we determine later. The subgradient of the objective function in the dual problem $(D)$ can be written in the following form 
\begin{equation}
\label{5_eq_grad}
g(p_1, \, \ldots, \, p_n)=\mathtt{x}(\p)-C\lambda(\p),
\end{equation}
where~$\mathtt{x}(\p)=(x_1(p_1), \, \ldots, \, x_n(p_n))^{\top}, \, \lambda(\p)=(\lambda_1(p_1), \, \ldots, \, \lambda_n(p_n))^{\top}$ and $\sum \limits_{k=1}^{n} \lambda_k(p_k)=1, \; \lambda_k(p_k) \geqslant 0$ if $k \, \in \, \argminset\limits_{j=1, \, \ldots, \, n} p_j$ and $\lambda_k(p_k) = 0$, if 
$k \, \notin \, \argminset\limits_{j=1, \, \ldots, \, n} p_j$.
Note that here and below $g (p_1, \, \ldots, \, p_n) \in \partial \varphi(p_1,\, \ldots, \,p_n) $ is an arbitrary subgradient, i.e. an arbitrary element of the convex compact set~--- subdifferential.

\begin{equation*}
\begin{array}{|c|}
\hline \\
\mbox{\bf General projected subgradient method}\\
\\
\hline \\
\quad \mbox{
\begin{minipage}{14cm}
\textbf{Input:} $h$ -- stepsize, $\p^0$ -- starting point.
\begin{enumerate}
	\item  Set $p_k^t, \, k=1, \, \ldots, \, n$ and calculate 
\begin{equation*}
x_k(p^{t}_k)=\argmax \limits_{x_k \geqslant  0} \Big \{ p^{t}_k x_k-f_k(x_k) \Big \}, \quad  k =1, \, 2, \, \ldots, \, n.
\end{equation*}
	\item 	Form a vector $\lambda(\p^t)$ as $\lambda^t=(\lambda^t_1, \, \ldots, \, \lambda^t_n)^{\top}$, where $\sum \limits_{k=1}^{n} \lambda_k^t=1$ and if $k \, \in \, \argminset\limits_{j=1, \, \ldots, \, n} p^t_j$, then $\lambda_k^t \geqslant 0$, otherwise, if $k \, \notin \, \argminset\limits_{j=1, \, \ldots, \, n} p_j^t$, then $\lambda_k^t = 0$
	\item  Do the step 
$$\p^{t+1}=
\left(\p^{t}-h g(p_1^t, \, \ldots, \,p_n^t) \right )_{+},$$
where  $$g(p_1^t, \, \ldots, \,p_n^t)=\mathtt{x}(\p^t) - C\lambda(\p^t).
$$
\end{enumerate}
\end{minipage}
}
\quad\\
\\
\hline
\end{array}
\end{equation*}

Next, we discuss the economic interpretation of the steps of the subgradient method.
First of all, $C\lambda_k(\p^{t})$ can be seen as the volume the Center intends to purchase from the $k$-th producer on $t$ iteration. Indeed, the dual variable $\p$ corresponds to the vector of prices set by the producers and non-zero components of the vector $\lambda$ can correspond only to those producers, who have minimal price. Since $\sum_{k=1: \; k \, \in \, \argminset\limits_{j=1, \, \ldots, \, n} p_j}^n \lambda_k^t = 1$ and all $\lambda_k^t$ are non-negative, $\lambda_k^t$ for $k\in\argminset\limits_{j=1, \, \ldots, \, n} p_j$ can be seen as the relative proportions of purchase form the producers having minimal price. Thus, $C\lambda_k(\p^{t})$ can be seen as the volume the Center intends to purchase from the $k$-th producer on $t$ iteration.

Further, each $k$-th component of the subgradient~\eqref{5_eq_grad} can be interpreted as the difference between the production $x_k(p_k)$ of the $k$-th factory and the volume $C\lambda_k$ of the Center's demand for this factory. 
For some factories, for which $C\lambda_k^t$ is positive, it can happen that $x_k(p_k)-C\lambda_k^t < 0$. This is a signal for the  $k$-th producer that the demand exceeds the supply and the $k$-th price can be increased together with the increase of the produced volume.

Finally, the subgradient step is the production adjustment steps for each producer, i.e. each producer counts how much its production differs from the desired volume of the Center's purchase from this producer this year. If the Center does not want to buy anything from the producer or buys less than it produced, then the producer lowers the price. If the Center is ready to buy more than the producer produced, the producer raises the price. In the case of equality, the producer does not change anything.

As a result the policy for the Center and producers is as follows.
\begin{equation}
\label{Alg1}
\begin{array}{|c|}
\hline \\
\mbox{\bf Subgradient method for the resource allocation}\\
\\
\hline \\
\quad \mbox{
\begin{minipage}{14cm}
\textbf{Input:} $\e > 0$ -- accuracy,  $\p^0$ -- starting point.
\begin{enumerate}
	\item Set the stepsize $h = \frac{\e}{nC^2}$. 
	\item  Given the price vector $\p^t$ for the current year, 
	producers calculate the optimal production plan for these prices as 
	\begin{equation*}
x_k(p^{t}_k)=\argmax \limits_{x_k \geqslant  0} \Big \{ p^{t}_k x_k-f_k(x_k) \Big \}, \quad  k =1, \, 2, \, \ldots, \, n.
\end{equation*}
and communicates this information to the Center.
	\item The Center determines the shares of purchases for each producer, i.e. forms a vector $\lambda(\p^t)$ as $\lambda^t=(\lambda^t_1, \, \ldots, \, \lambda^t_n)^{\top}$, where $\sum \limits_{k=1}^{n} \lambda_k^t=1, \; \lambda_k^t \geqslant 0$ if $k \, \in \, \argminset\limits_{j=1, \, \ldots, \, n} p^t_j$ and $\lambda_k^t = 0$, if 
$k \, \notin \, \argminset\limits_{j=1, \, \ldots, \, n} p_j^t$
	and sends this vector to all factories.
	
	\item  Each factory adjusts the price for the next year as follows
$$\p^{t+1}=
\left(\p^{t}-h (\mathtt{x}(\p^t) - C\lambda(\p^t)) \right )_{+}.$$

\end{enumerate}
\end{minipage}
}
\quad\\
\\
\hline
\end{array}
\end{equation}

To state the convergence rate result, we need introduce an upper bound for the optimal value of the prices. 

\begin{Lm}
\label{Lm:1}
Let the $\p^{*}$ be a solution to the dual problem $(D).$ Then 
\begin{equation}
\label{estim_p^*}
   \norm{\p^{*}}_{2}\leqslant \sqrt{n}p_{max}.
\end{equation}
where 
\begin{equation}
  \label{eq:p_max}
   p_{max} := \dfrac{n}{C}\left(\sum \limits_{k=1}^n f_k\left(\dfrac{2C}{n}\right)-\sum \limits_{k=1}^n f_k(0)\right). 
\end{equation}
\end{Lm}
The proof of this lemma is deferred to the Appendix.

Then, we can formulate the following theorem about convergence rate
\begin{Th}[\cite{Ivanova}]
Let Algorithm \eqref{Alg1} be run with starting point $\mathtt{p}^0$ satisfying $0 \leqslant p_k^0 \leqslant p_{max},\; k=1, \, \ldots, \, n$ for 
$$
N= \left\lceil\dfrac{164(Cnp_{max})^2}{\varepsilon^2}\right\rceil 
$$
steps.
Then
\begin{equation}
\label{stop_crit}
f(\mathtt{x}^N)-f(\mathtt{x}^{*})\leqslant   \varepsilon,\;  C -\sum\limits_{k=1}^{n}x_{k}^N\leqslant   \dfrac{\varepsilon}{3p_{max}},
\end{equation}
where $\mathtt{x}^N = \dfrac{1}{N}\sum\limits_{t=0}^{N-1}\mathtt{x}(\p^{t})$
\end{Th}
Note that the number of iterations $N$ to achieve accuracy $\e$ is very large. 
To improve the number of iterations, in the following sections we consider the methods based on the composite optimization approach.

\section{Composite gradient method for the resource allocation problem}
\label{comp}
In this section we consider a non-accelerated composite gradient method to solve the dual problem $(D)$, including its interpretation and convergence rate estimate. The problem $(D)$ can be rewritten as
\begin{equation*}
\varphi(p_1, \, \ldots, \, p_n)=\psi(p_1, \, \ldots, \, p_n)+g(p_1, \, \ldots, \, p_n),
\end{equation*}
where
\begin{equation}
\label{tildephi}
\psi(p_1, \, \ldots, \, p_n)=\sum \limits_{k=1}^n \Big\{ p_{k}x_k(p_k)-f_k(x_k(p_k))\Big\}=\left<\p,\mathtt{x}(\p)\right>-f(\mathtt{x}(\p))
\end{equation}
is convex function. Gradient of this function
\begin{equation}
\label{psi_grad}
\nabla \psi(p)=\mathtt{x}(\p),
\end{equation}
satisfies Lipschitz condition (see e.g. \cite{Smooth})
\begin{equation}
\label{Lips_grad}
\norm{\nabla\psi(\p^{1})-\nabla\psi(\p^{2})}_2\leqslant L_{\psi}\norm{\p^{1}-\p^{2}}_2, \, \, \forall \, \p^{1}, \p^{2} \geqslant 0,
\end{equation}
where $L_{\psi}=\dfrac{n}{\mu}$ (see e.g.\cite{Ivanova}) and 
\begin{equation*}
g(p_1, \, \ldots, \, p_n)=-C\min \limits_{k=1, \, \ldots, \, n}p_k
\end{equation*}
is convex non smooth function. 
The idea is to use first-order information about $\psi$ and use the function $g$ as a whole, see the main projected gradient step \eqref{step_1}. This allows the method to work in accordance with the smoothness properties of $\psi$, i.e. be faster since $\psi$ is smooth. The price is that the step \eqref{step_1} has to be simple enough. The next lemma shows that this is indeed the case and the solution can be found explicitly.

\begin{equation*}
\begin{array}{|c|}
\hline \\
\mbox{\bf General composite projected gradient method}\\
\\
\hline \\
\quad \mbox{
\begin{minipage}{14cm}
\textbf{Input:} $N > 0$ -- number of steps, $L_{\psi}$ -- Lipschitz constant of gradient $\psi$, $\p^0$ -- starting point.
\begin{enumerate} 
	\item Find 
\begin{equation*}
x_k(p^{t}_k)=\argmax \limits_{x_k \geqslant  0} \Big \{ p^{t}_k x_k-f_k(x_k) \Big \}, \quad  k =1, \, 2, \, \ldots, \, n.
\end{equation*}
	
	\item  Do the step 
\begin{equation}
\label{step_1}
    \p^{t+1}=\argmin_{\p\geqslant 0}\left\{\left<\nabla\psi(\p^t),\p-\p^{t}\right>-C\min \limits_{k=1, \, \ldots, \, n}p_k+\dfrac{L_{\psi}}{2}\norm{\p-\p^t}_{2}^{2}\right\}.
\end{equation}
\end{enumerate}
\end{minipage}
}
\quad\\
\\
\hline
\end{array}
\end{equation*}

\begin{Lm}
\label{lemma2}
Let $\tilde{\p}^{t+1}=\p^t-\dfrac{1}{L_{\psi}}\mathtt{x}(\p^t)$. Then ${\p}^{t+1}$ in~\eqref{step_1} is defined as follows
    $$
    p^{t+1}_{k}=\max\left(p^{t+1}_{center},\,  \tilde{p}^{t+1}_k\right), \, \, k=1,\, \ldots\, ,n,
    $$
where if $\sum\limits_{k=1}^{n}\left(-\tilde{p}^{t+1}_k \right)_{+}\geqslant \dfrac{C}{ L_{\psi}}$ than $p_{center}^{t+1}=0$,  else $p^{t+1}_{center} > 0$ is a solution of equation
    $$
    \sum\limits_{k=1}^{n}\left(p^{t+1}_{center}- \tilde{p}^{t+1}_k \right)_{+}= \dfrac{C}{ L_{\psi}}.
    $$

\end{Lm}
The proof of Lemma~\ref{lemma2} can be found in Appendix.
Note that the step~\eqref{step_1} can be rewritten as 
\begin{equation}
\label{step}
    \p^{t+1}=\left[\p^{t}-\dfrac{1}{L_{\psi}}\left(\mathtt{x}(\p^t)-C\lambda(\p^{t+1})\right)\right]_{+},
\end{equation}
where $\lambda(\p^{t+1})$ is such that 
$\sum \limits_{k=1}^{n} \lambda_k(p_k^{t+1})=1, \; \lambda_k(p_k^{t+1}) \geqslant 0$ if $k \, \in \, \argminset\limits_{j=1, \, \ldots, \, n} p_j^{t+1}$ and $\lambda_k(p_k^{t+1}) = 0$, if 
$k \, \notin \, \argminset\limits_{j=1, \, \ldots, \, n} p_j^{t+1}$.
This equality looks very similar to equality~\eqref{5_eq_subgrad_der_1}.

\subsection{Composite gradient method for the resource allocation problem}
In this subsection, we apply general projected gradient method to the resource allocation problem and give its interpretation, i.e. describe the strategy of the Center and producers, and state the convergence rate theorem for this method. 
First of all, from \eqref{Lips_grad}, since $\nabla\psi(\p)=\mathtt{x}(\p)$, the constant $L_{\psi}$ determines the relation between prices $\p$ and production $\mathtt{x}(\p)$, i.e. says, what is the maximum change in the production if the price changes. 
From Lemma~\ref{lemma2}, each year $t$ the Center, knowing the prices and production of each producer, forms a prediction $\tilde{p}^{t+1}_k=p^t_k-\dfrac{1}{L_{\psi}}x_k(p_k^t)$ for the lowest possible producers' prices vector for the next year. 
After that, the goal of the Center is to set its purchase price $p^{t+1}_{center}$ to satisfy the total demand $C$. To explain, how it is done, let us look at \eqref{step}, which looks very similar to~\eqref{5_eq_subgrad_der_1}. The key difference is that in \eqref{step} $\lambda$ depends on the unknown price vector $p_k^{t+1}$, resulting in the implicit definition of this price vector as a solution to nonlinear equation. 
But, $C\lambda(\p^{t+1})$ is still an estimate of the purchase from each producer for the next year and the Center chooses to buy only from the producers having the lowest price. By Lemma~\ref{lemma2}, the solution of~\eqref{step} can be written as   
\begin{equation*}
    \p^{t+1}=\left[\tilde{\p}^{t+1}+\dfrac{C}{L_{\psi}}\lambda(\p^{t+1})\right]_{+},
\end{equation*}
where 
$$
\lambda_k(p_k^{t+1}) = \dfrac{L_{\psi}}{C}\left(p^{t+1}_{center}- \tilde{p}^{t+1}_k \right)_{+}.
$$
Thus, the estimate for the purchase from the $k$-th producer is $L_{\psi}\left(p^{t+1}_{center}- \tilde{p}^{t+1}_k \right)_+$, i.e. if the predicted price $\tilde{p}^{t+1}_k$ is higher than $p^{t+1}_{center}$, there will be no purchase from this producer. Again, by Lemma~\ref{lemma2}, if for the smallest possible price $p^{t+1}_{center} = 0$, the total purchased amount $L_{\psi}\sum\limits_{k=1}^{n}\left(-\tilde{p}^{t+1}_k \right)_{+}$ is greater than $C$, the Center sets its price $p^{t+1}_{center} = 0$. Otherwise, it determines such price $p^{t+1}_{center}$ that the total estimated purchase $L_{\psi}\sum\limits_{k=1}^{n}\left(p^{t+1}_{center}-\tilde{p}^{t+1}_k \right)_{+}$ is exactly $C$. Having defined $p^{t+1}_{center}$, the Center informs the producers about this price.

From the perspective of the producer, similarly to section~\ref{Subgradient descent}, each producer $k$, knowing the price $p_k^{t}$ for the current year $t$, determines the optimal production plan $x_k(p^t_k)$ for this price and reports the price and production to the Center. 
After receiving $p^{t+1}_{center}$ from the center, each producer adjusts its price $p_k^{t+1}$ for the next year as  $p^{t+1}_{k}=\max\left(p^{t+1}_{center},\,  \tilde{p}^{t+1}_k\right)$. This means that the lowest possible price $\tilde{p}^{t+1}_k$ is compared to the Center's price. If $p^{t+1}_{center} \geqslant \tilde{p}^{t+1}_k$, the producer can increase the price up to $p^{t+1}_{center}$ and the Center will buy. On the other hand, if $p^{t+1}_{center} < \tilde{p}^{t+1}_k$, there is no sense for the producer to sell at the price $p^{t+1}_{center}$ and have a loss.

We summarize the strategies of the Center and producers in Algorithm~\eqref{Alg2}.

Note that at equilibrium the Center will purchase from all producers and its price will be constant at all iterations. That is, all producers will be able to set a minimum price for the next year so that the Center buys from them the optimal volume.

\begin{equation}
\label{Alg2}
\begin{array}{|c|}
\hline \\
\mbox{\bf Composite gradient method for the resource allocation}\\
\\
\hline \\
\quad \mbox{
\begin{minipage}{13cm}
\textbf{Input:} $N> 0$ -- number of steps, $L_{\psi}$ -- Lipschitz constant of gradient $\psi$, $\p^0$ -- starting point.
\begin{enumerate}
	\item  Knowing the prices $p_k^t, \, k=1, \, \ldots, \, n$ for the current year $t$, producers calculate  the optimal plan for the production according these prices 
	\begin{equation*}
x_k(p^{t}_k)=\argmax \limits_{x_k \geqslant  0} \Big \{ p^{t}_k x_k-f_k(x_k) \Big \}, \quad  k =1, \, 2, \, \ldots, \, n.
\end{equation*}
\item The Center forms a prediction for the lowest possible producers' prices vector for the next year 
$$
\tilde{p}^{t+1}_k=p^t_k-\dfrac{1}{L_{\psi}}x_k(p_k^t), \quad  k =1, \, 2, \, \ldots, \, n.
$$
	\item The Center determines the price $p^{t+1}_{center}$ at which it will purchase product for the next year $t+1$ and sends this price to all factories.  
\begin{itemize}
    \item If $\sum\limits_{k=1}^{n}\left(-\tilde{p}^{t+1}_k \right)_{+}\geqslant \dfrac{C}{ L_{\psi}}$ then $p^{t+1}_{center}=0$;
    \item Else $p^{t+1}_{center} > 0$ and solves equation 
    $$
    \sum\limits_{k=1}^{n}\left(p^{t+1}_{center}-\tilde{p}^{t+1}_k \right)_{+}= \dfrac{C}{ L_{\psi}}
    $$
\end{itemize}

	\item  Each producer adjusts the price for the next year as follows
$$
    p^{t+1}_{k}=\max\left(p^{t+1}_{center},\, \tilde{p}^{t+1}_k\right), \, \, k=1,\, \ldots\, ,n.
    $$

\end{enumerate}
\end{minipage}
}
\quad\\
\\
\hline
\end{array}
\end{equation}

\begin{Th}
\label{th2}
Let Algorithm \eqref{Alg2} be run for $N$ steps with starting point $\mathtt{p}^0$ satisfying $0 \leqslant p_k^0 \leqslant p_{max},\; k=1, \, \ldots, \, n$, where $p_{max}$ is given in \eqref{eq:p_max}.
Then
\begin{eqnarray*}
f(\mathtt{x}^N)-f(\mathtt{x}^{*})&\leqslant   f(\mathtt{x}^N)+\varphi(\mathtt{p}^{*})\leqslant   \varphi(\mathtt{p}^N)+f(\mathtt{x}^N) \leqslant \dfrac{82p_{max}^2n^2}{N\mu} , \notag\\
 \left[C- \sum\limits_{k=1}^{n}x_k^N \right]_{+} & \leqslant \dfrac{82p_{max}n^2}{3N\mu},\notag
\end{eqnarray*}
where $\p^{N}=\dfrac{1}{N}\sum\limits_{t=1}^{N}\p^{t}$ and $\mathtt{x}^{N}=\dfrac{1}{N}\sum\limits_{t=0}^{N-1}\mathtt{x}(\p^{t})$.
\end{Th}
The proof of this theorem is deferred to the Appendix.

Let us make a remark on the complexity of this procedure. Assume that we want to solve problem $(P)$ with accuracy $\e$ in the following sense
\begin{equation}
\label{stop_crit_1}
f(\mathtt{x}^N)-f(\mathtt{x}^{*})\leqslant   \varepsilon,\;  C -\sum\limits_{k=1}^{n}x_{k}^N\leqslant   \dfrac{\varepsilon}{3p_{max}}.
\end{equation}
To solve this problem, we consider the dual problem $(D)$, and solve it by the composite gradient method starting with  $\mathtt{p}^0$ satisfying $0 \leqslant p_k^0 \leqslant p_{max},\; k=1, \, \ldots, \, n$, where $p_{max}$ is given in \eqref{eq:p_max}.

Theorem \ref{th2} states that Algorithm will find the solution 
no later than after $N=\dfrac{82p_{max}^2n^2}{\varepsilon\mu}$ iterations.

\section{Accelerated composite gradient method for the resource allocation problem}
\label{AGM}
In this section we use accelerated composite gradient method to solve the resource allocation problem. Accelerated algorithm allows to improve the convergence rate in comparison to previous section. 

We can rewrite \ref{eq:PDzetakp1Def} step as follows 
\begin{eqnarray*}
\mathtt{\y}^{t+1}&=&\argmin_{\p\geqslant 0}\left\{\alpha_{t+1}\left(\left<\nabla\psi(\p^{t+1}),\p-\p^{t+1}\right>-C\min \limits_{k=1, \, \ldots, \, n}p_k\right)+\dfrac{1}{2}\norm{\p-\y^t}_{2}^{2}\right\}\\&
=&\argmin_{\p\geqslant 0}\left\{\ \left<\mathtt{x}(\p^{t+1}),\p-\p^{t+1}\right>+C\max \limits_{k=1, \, \ldots, \, n}(-p_k)+\dfrac{1}{2\alpha_{t+1}}\norm{\p-\y^t}_{2}^{2}\right\}\\&
=&\argmin_{\p\geqslant 0}\left\{\ \max \limits_{k=1, \, \ldots, \, n}(-p_k)+\dfrac{1}{2C\alpha_{t+1}} \norm{\p-\left(\y^t-\alpha_{t+1}\mathtt{x}(\p^{t+1})\right)}_{2}^{2}\right\}.
\end{eqnarray*}
Define $\tilde{\y}^{t+1}=\y^t-\alpha_{t+1}\mathtt{x}(\p^{t+1})$ 
and then, using Lemma~\ref{lemma2}, we obtain the following solution:
\begin{itemize}
    \item If $\sum\limits_{k=1}^{n}\left(-\tilde{y}^{t+1}_k \right)_{+} \geqslant C\alpha_{t+1}$ then $y^{t+1}_{center}=0$ and 
    $$
    y^{t+1}_{k}=\max\left(0,\, \tilde{y}^{t+1}_k\right), \, \, k=1,\, \ldots\, ,n.
    $$
    \item Else $y^{t+1}_{center} > 0$ is determined from 
    $$
    \sum\limits_{k=1}^{n}\left(y^{t+1}_{center}- \tilde{y}^{t+1}_k \right)_{+}= C\alpha_{t+1}
    $$
    and 
    $$
    y^{t+1}_{k}=\max\left(y^{t+1}_{center},\,  \tilde{y}^{t+1}_k\right), \, \, k=1,\, \ldots\, ,n.
    $$
\end{itemize}
\begin{equation*}
\begin{array}{|c|}
\hline \\
\mbox{\bf Accelerated composite gradient descent scheme}\\
\\
\hline \\
\quad 
\mbox{
\begin{minipage}{15cm}
\textbf{Input:} $N> 0$ -- number of steps, $L_{\psi}$ -- Lipschitz constant of gradient $\psi$, $\p^0, \, \y^0, \, \w^0$ -- starting points, $\alpha_0=A_0=0$.
\begin{enumerate}
	\item  Find $\alpha_{t+1}$ as the largest root of the equation
				\begin{equation}
				\label{PDalpQuadEq}
				A_{t+1}:=A_t+\alpha_{t+1} = L_{\psi}\alpha_{t+1}^2.
				\end{equation}
\item Calculate 
				\begin{equation}
				\p^{t+1} = \frac{\alpha_{t+1}\y^t + A_t \w^t}{A_{t+1}} .
				\label{eq:PDDef1}
				\end{equation}
	\item Calculate 
				\begin{equation}
			\mathtt{\y}^{t+1}=\argmin_{\p\geqslant 0}\left\{\alpha_{t+1}\left(\left<\nabla\psi(\p^{t+1}),\p-\p^{t+1}\right>+g(\p)\right)+\dfrac{1}{2}\norm{\p-\y^t}_{2}^{2}\right\}.
				\label{eq:PDzetakp1Def}
				\end{equation}
				
	\item Calculate
				\begin{equation}
				\w^{t+1} = \frac{\alpha_{t+1}\y^{t+1} + A_t \w^t}{A_{t+1}}.
				\label{eq:PDDef2}
				\end{equation}

	\item Set 
				\begin{equation*}
					\mathtt{x}^{t+1} = \frac{1}{A_{t+1}}\sum_{i=0}^{t+1} \alpha_i \mathtt{x}(\p^i) = \frac{\alpha_{t+1}x(\p^{t+1})+A_t\mathtt{x}^{t}}{A_{t+1}}.
				\notag
				\end{equation*}
\end{enumerate}
\end{minipage}
}
\quad\\
\\
\hline
\end{array}
\end{equation*}

We apply this general method to the resource allocation problem. Let us make some remarks on the interpretation.
	
Note that in the non-accelerated method we have fixed $L_{\psi}$, which we interpreted as the expectation of the Center on how the amount of production (of any producer) changes, if the purchase price is changed by one. This was the upper bound for this value. In the accelerated method this role is played by $\alpha_t$, which is adjusted on each iteration and allows to have faster convergence.
Vector $\w^{t+1}$ can be interpreted as an average of the historical price data $\y^t$, $t\geq 0$. As we prove in next theorem, $\w^{t+1}$ is a good approximation for the solution of the dual problem. 
This price can give a hint to each producer on the value of the predicted price for the next period. Thus, this average historical price is combined with the prediction price $\y^t$ from the previous period to form $\p^{t+1}$. Then the optimal production plan is calculated and, similarly to the non-accelerated method, producers predict the price for next year $\tilde{\y}^{t+1}$ and communicate it to the Center. After getting from the Center $y^{t+1}_{center}$, each producer adjust the prediction price for the next year.

\begin{Th}
\label{th3}
Let Algorithm \eqref{Alg3} be run for $N$ steps with starting points $\p^0=\y^0=\w^0$ satisfying $0 \leqslant p_k^0 \leqslant p_{max},\; k=1, \, \ldots, \, n$, where $p_{max}$ is given in \eqref{eq:p_max}. Then
\begin{align}
f(\mathtt{x}^N)-f(\mathtt{x}^{*})&\leqslant   f(\mathtt{x}^N)+\varphi(\mathtt{w}^{*})\leqslant   \varphi(\mathtt{w}^N)+f(\mathtt{x}^N) \leqslant \dfrac{148n^2p_{max}^2}{(N+1)^2\mu} , \notag \\
 \left[C- \sum\limits_{k=1}^{n}x_k^N \right]_{+} & \leqslant \dfrac{148n^2Rp_{max}}{5(N+1)^2\mu}, \notag
\end{align}
where $\mathtt{x}^N=\dfrac{1}{A_{N}}\sum\limits_{t=0}^{N}\alpha_t\mathtt{x}(\p^t)$
\end{Th}
The proof of this theorem mostly follows the steps of the Theorem~\ref{th2} in the previous section, but we give all steps for the reader’s convenience in Appendix.

\section{Numerical experiments}
In this section we compare and verify the algorithms presented in the paper. All the algorithms were implemented in the programming language~Python. 

\begin{figure}[H]
    \includegraphics[width=0.4\textwidth]{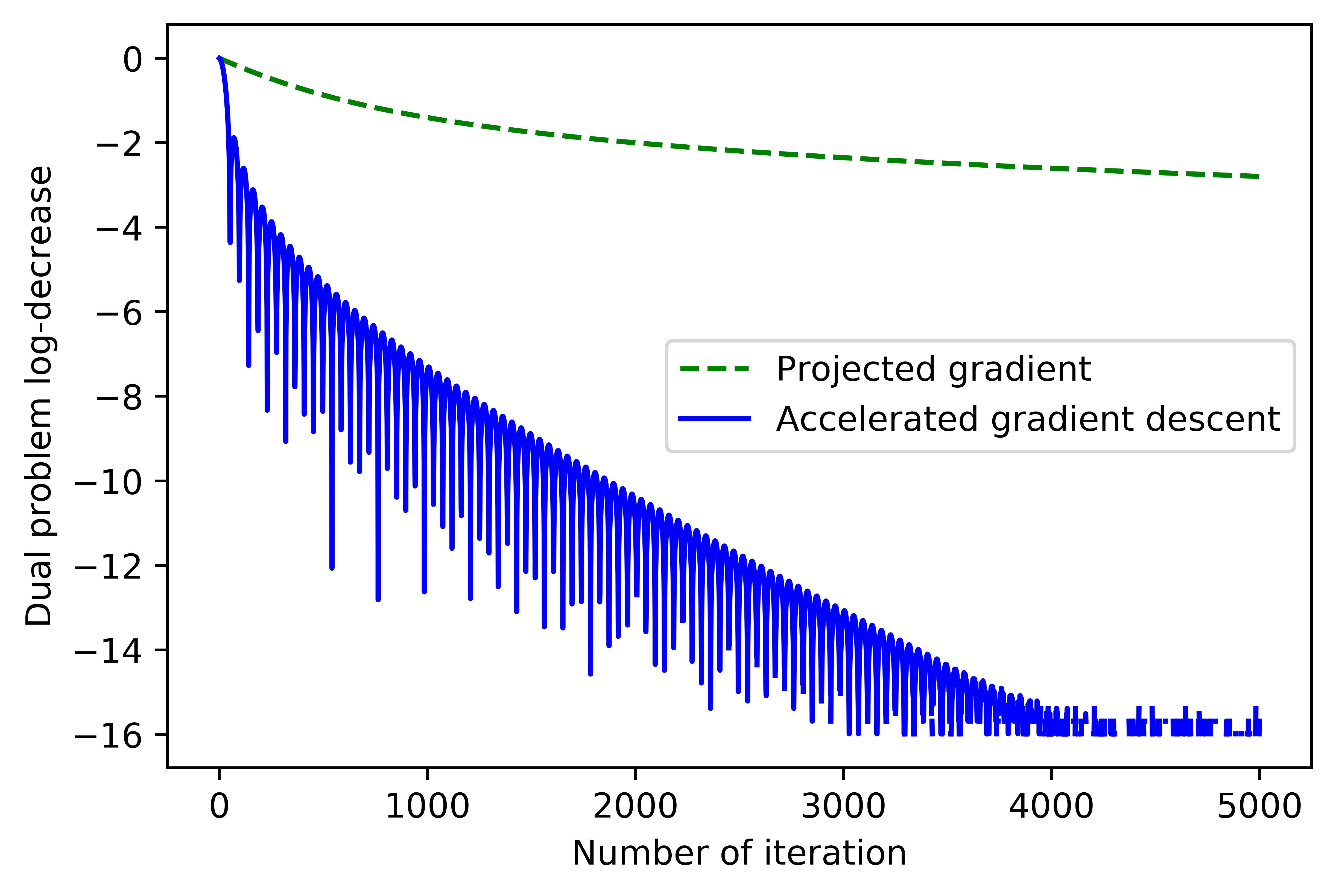}
    \includegraphics[width=0.4\textwidth]{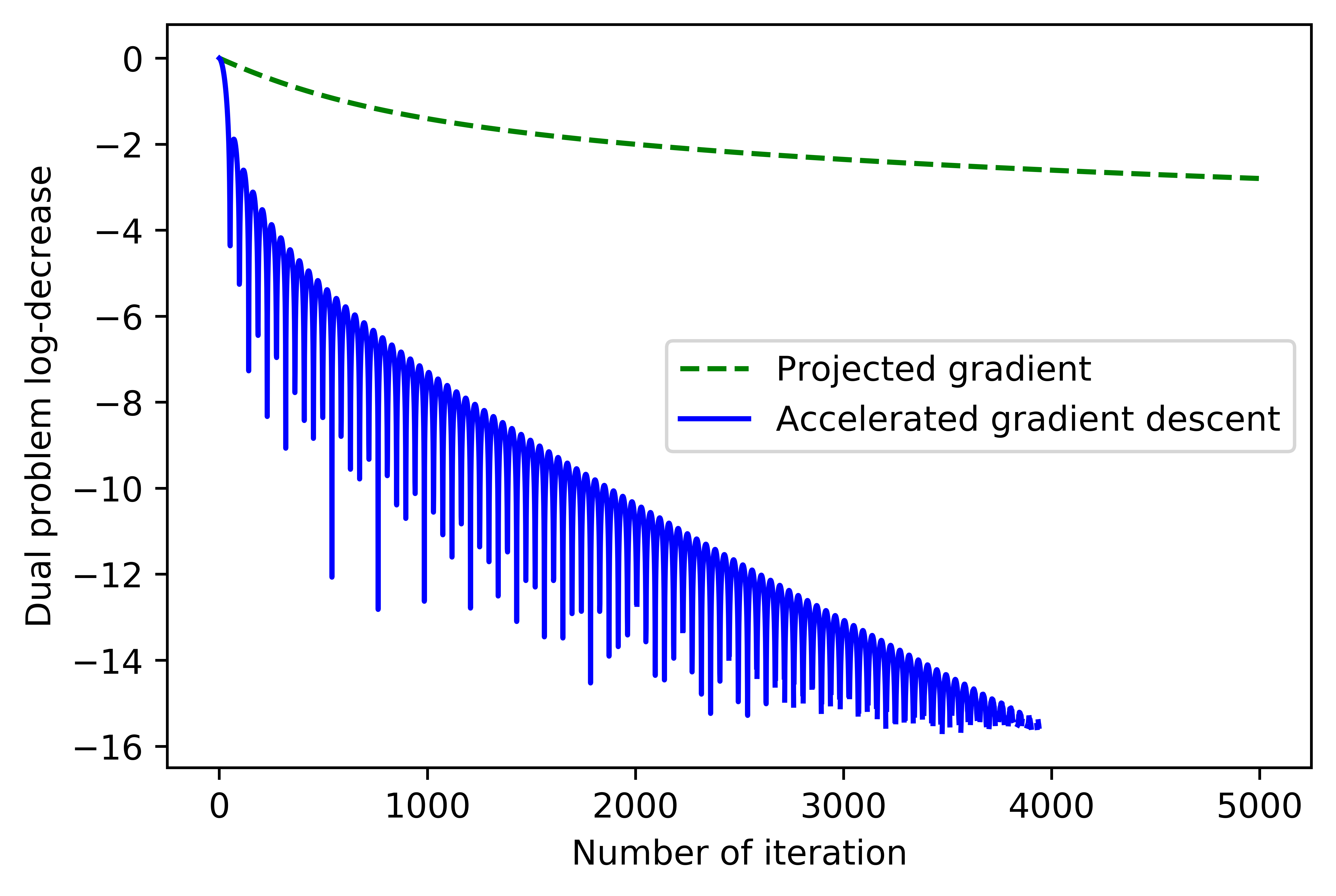}
\caption{Projected gradient vs Accelerated gradient descent,
 $100$~factories, cost functions~\eqref{eq_p5_prob1_fc1}, Dual value in log scale, 1 run and mean of 20 run}
\label{dr_graph1}
\end{figure}

\begin{figure}[H]              \includegraphics[width=0.4\textwidth]{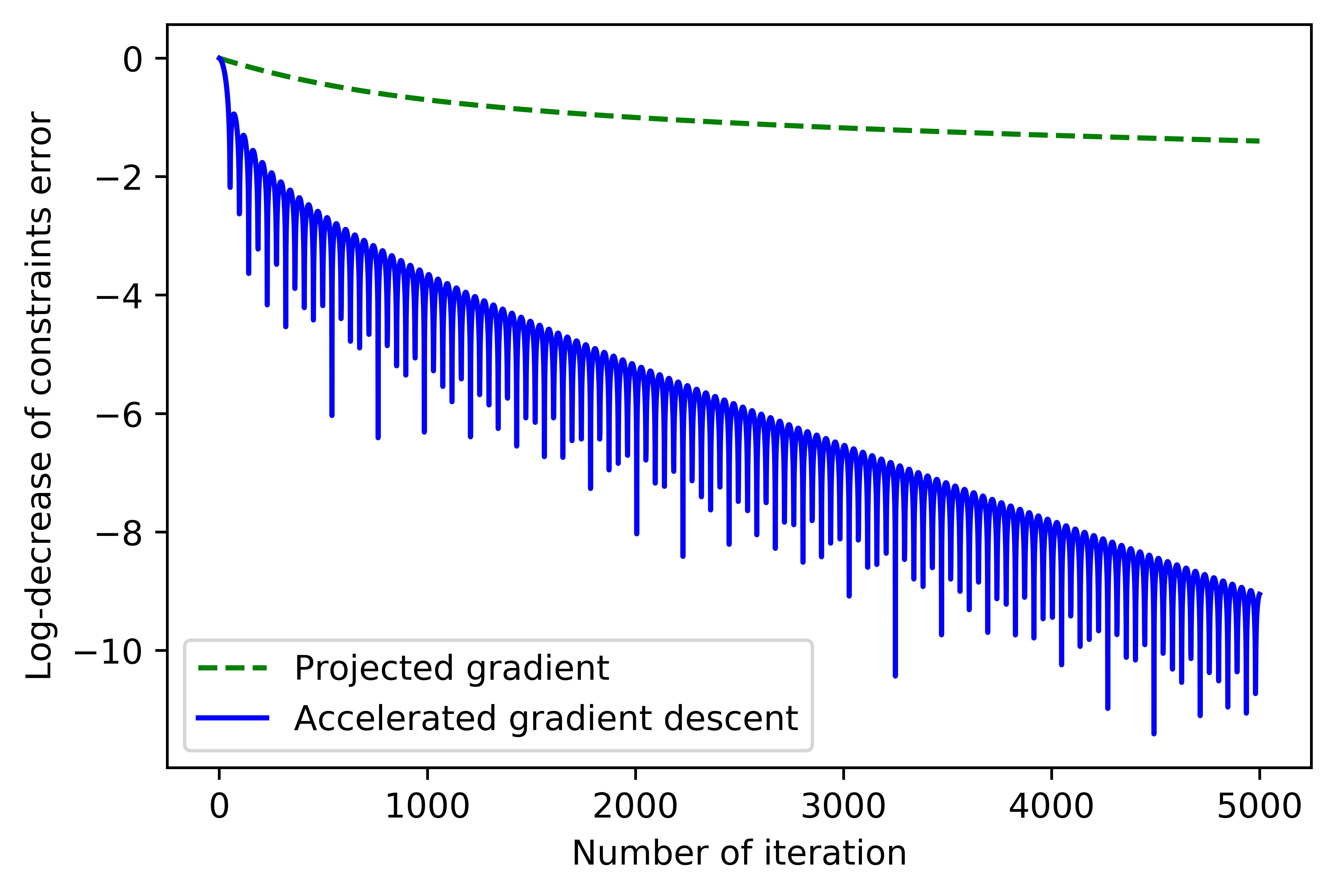}
    \includegraphics[width=0.4\textwidth]{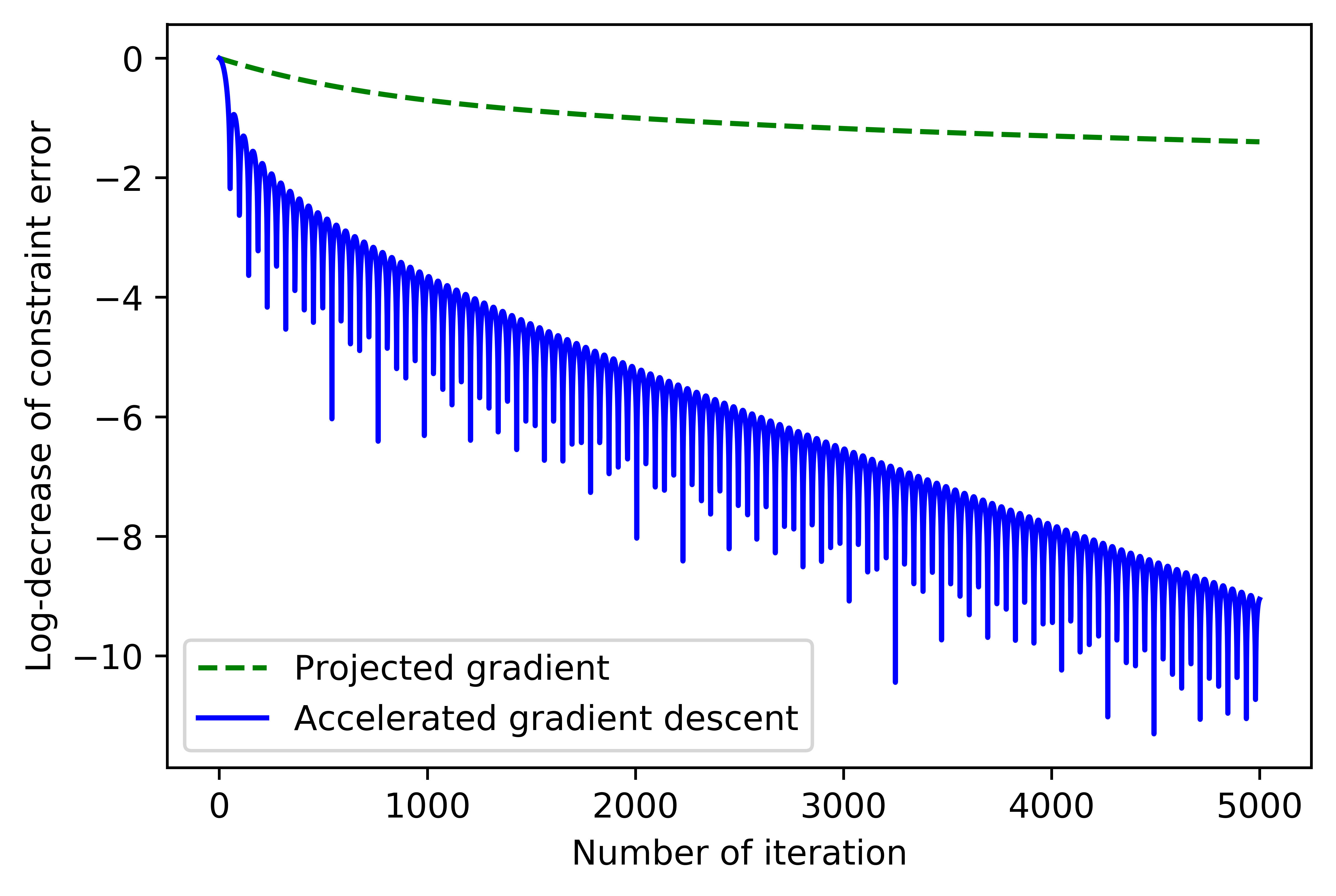}
\caption{Projected gradient vs Accelerated gradient descent,
 $100$~factories, cost functions~\eqref{eq_p5_prob1_fc1}, Inexactness of constraints in log scale, 1 run and mean of 20 run}
\label{dr_graph2}
\end{figure}

\begin{figure}[H]
\includegraphics[width=0.4\textwidth]{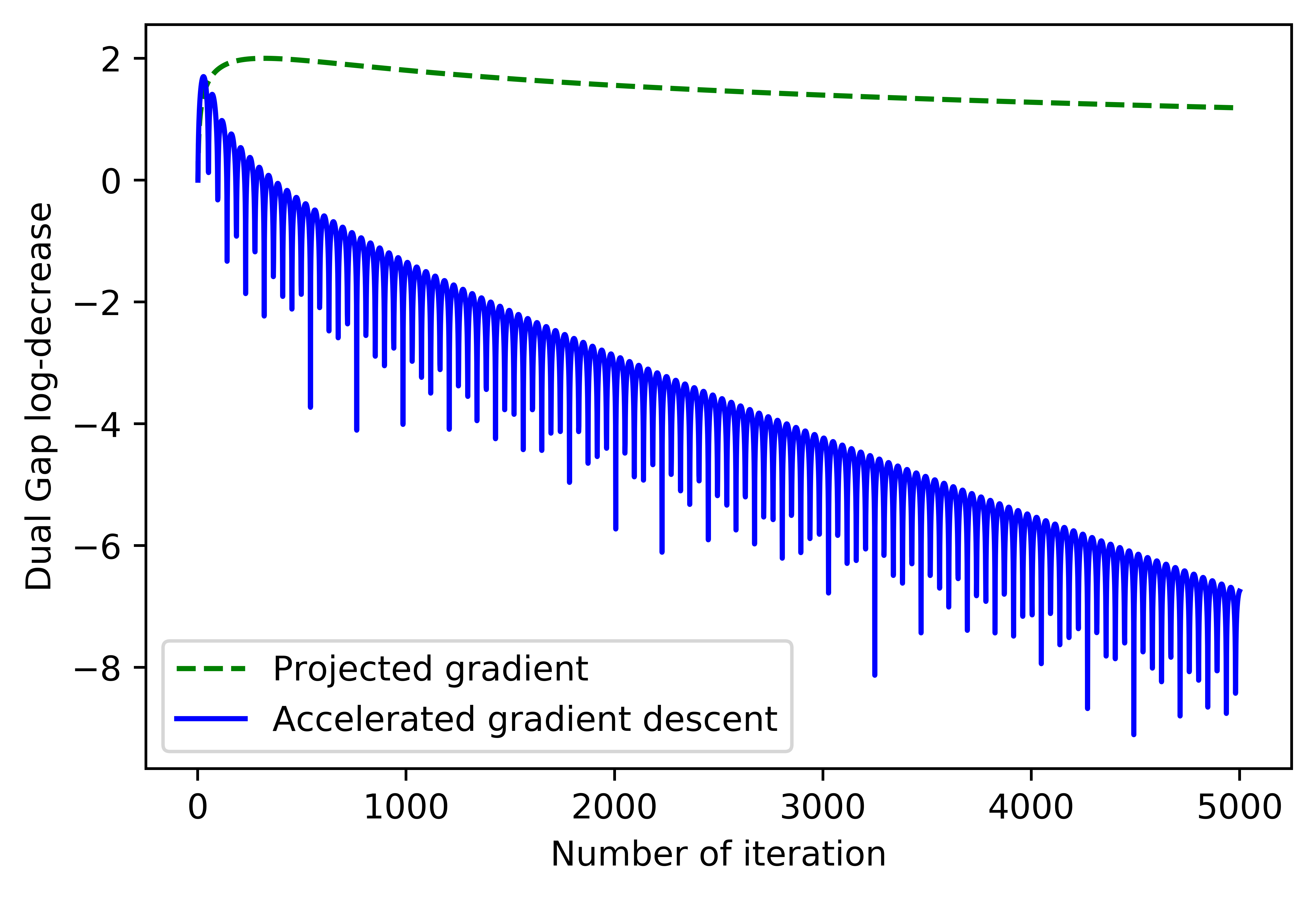}
    \includegraphics[width=0.4\textwidth]{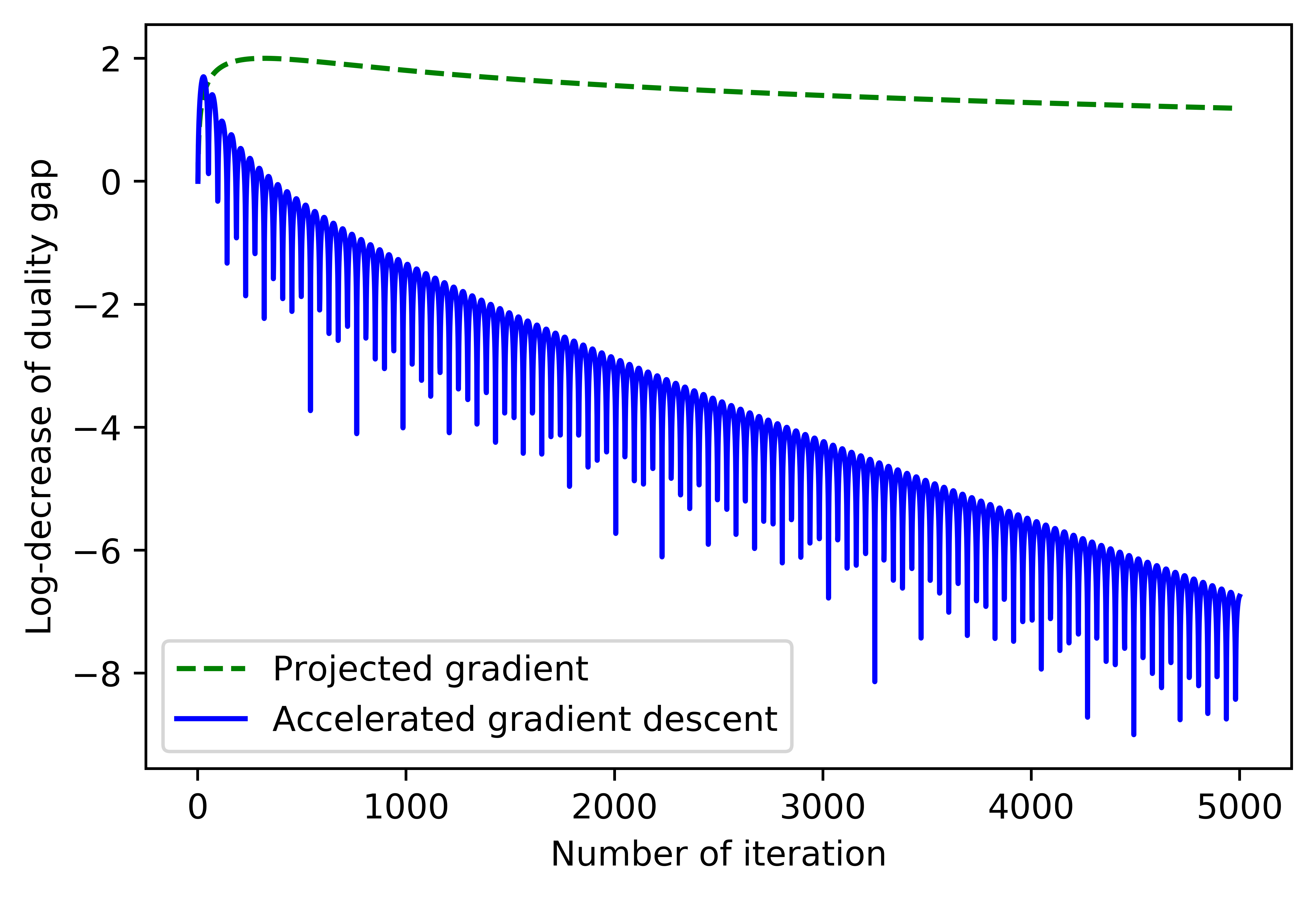}
\caption{Projected gradient vs Accelerated gradient descent,
 $100$~factories, cost functions~\eqref{eq_p5_prob1_fc1}, Value of duality gap in log scale, 1 run and mean of 20 run}
\label{dr_graph3}
\end{figure}

Firstly, we compare composite projected gradient method and accelerated composite gradient descent on a simple problem.
We assume that there are $100$ local producers of some long-living product, for example, wood. All this producers sell it to pulp and paper company (Center). Center needs to buy ~$10 000$ tons of wood. Each local producer has the following cost function: 
\begin{equation}
\label{eq_p5_prob1_fc1}
f_k(x_k) = \alpha_k x_k+\frac{\mu}{2} x_k^2, \quad k = 1, \, \ldots, \, 100,
\end{equation}
where $\alpha_k$ is the different local costs of production and transportation for each producer, $\mu$ is the costs of production that depends on technology and risks, that's why the same for producers. We consider this problem with random $\alpha_k \in N[100;400]$, $\mu=2$. Next in graphics the left side will be experiments for one start with random $\alpha$ and the right side will be the mean of $20$ random starts. The practical convergence rate for the objective value of the dual problem (D) is shown in log scale in fig.~\ref{dr_graph1}. The convergence rate for the constraint feasibility is shown in log scale in fig.~\ref{dr_graph2}.
The convergence rate for the value of the duality gap is shown in log scale  in fig.~\ref{dr_graph3}.

\begin{equation}
\label{Alg3}
\begin{array}{|c|}
\hline \\
\mbox{\bf Accelerated composite gradient descent for the resource allocation}\\
\\
\hline \\
\quad \mbox{
\begin{minipage}{15cm}
\textbf{Input:} $N> 0$ -- number of steps, $L_{\psi}$ Lipschitz constant of gradient $\psi$, $\p^0=\y^0=\w^0$ -- starting point.
\begin{enumerate}
	\item  In the current year $t$ producers find $\alpha_{t+1}$ as the largest root of the equation
				\begin{equation*}
				A_{t+1}:=A_t+\alpha_{t+1} = L_{\psi}\alpha_{t+1}^2.
				\end{equation*}

		\item All producers calculate the average price $\p^{t+1}$ as the following convex combination 
				\begin{equation*}
				p_k^{t+1} = \frac{\alpha_{t+1}y^t_k + A_t w_k^t}{A_{t+1}}, \, \, k=1,\, \ldots\, ,n 
				\end{equation*}
	 and calculate  the optimal plan for the production according to this price as 
	\begin{equation*}
x_k(p^{t+1}_k)=\argmax \limits_{x_k \geqslant  0} \Big \{ p^{t+1}_k x_k-f_k(x_k) \Big \}, \quad  k =1, \, 2, \, \ldots, \, n.
\end{equation*}
\item Each producer predicts the price for the next year $t+1$  as 
$$
\tilde{y}_k^{t+1}=y^t_k-\alpha_{t+1}x_k(\p^{t+1}), \, \, k=1,\, \ldots\, ,n
$$
and sends this information to the Center.
	\item The Center determines the prediction price $y^{t+1}_{center}$ at which it will purchase product for the next year $t+1$ as 
\begin{itemize}
    \item If $\sum\limits_{k=1}^{n}\left(-\tilde{y}^{t+1}_k \right)_{+} \geqslant C\alpha_{t+1}$ then $y^{t+1}_{center}=0$ 
    \item Else $y^{t+1}_{center} > 0$ and is determined from 
    $$
    \sum\limits_{k=1}^{n}\left(y^{t+1}_{center}- \tilde{y}^{t+1}_k \right)_{+}= C\alpha_{t+1},
    $$
    \end{itemize}
  
and sends this price to all producers.
	
	\item  Each producer adjusts the prediction price for the next year as 
 $$
    y^{t+1}_{k}=\max\left(y^{t+1}_{center},\,  \tilde{y}^{t+1}_k\right), \, \, k=1,\, \ldots\, ,n.
    $$
and calculates the historical price for for the next year 
				\begin{equation*}
				w_k^{t+1} = \frac{\alpha_{t+1}y_k^{t+1} + A_t w_k^t}{A_{t+1}}, \, \, k=1,\, \ldots\, ,n.
				\end{equation*}
\end{enumerate}
\end{minipage}
}
\quad\\
\\
\hline
\end{array}
\end{equation}

From the plots we can see that both methods work well. Projected gradient is monotone and converges in accordance with theoretical convergence rates. Accelerated gradient descent is non-monotone, faster than projected gradient and converges in accordance with theoretical convergence rates.

\begin{figure}[H]    \includegraphics[width=0.4\textwidth]{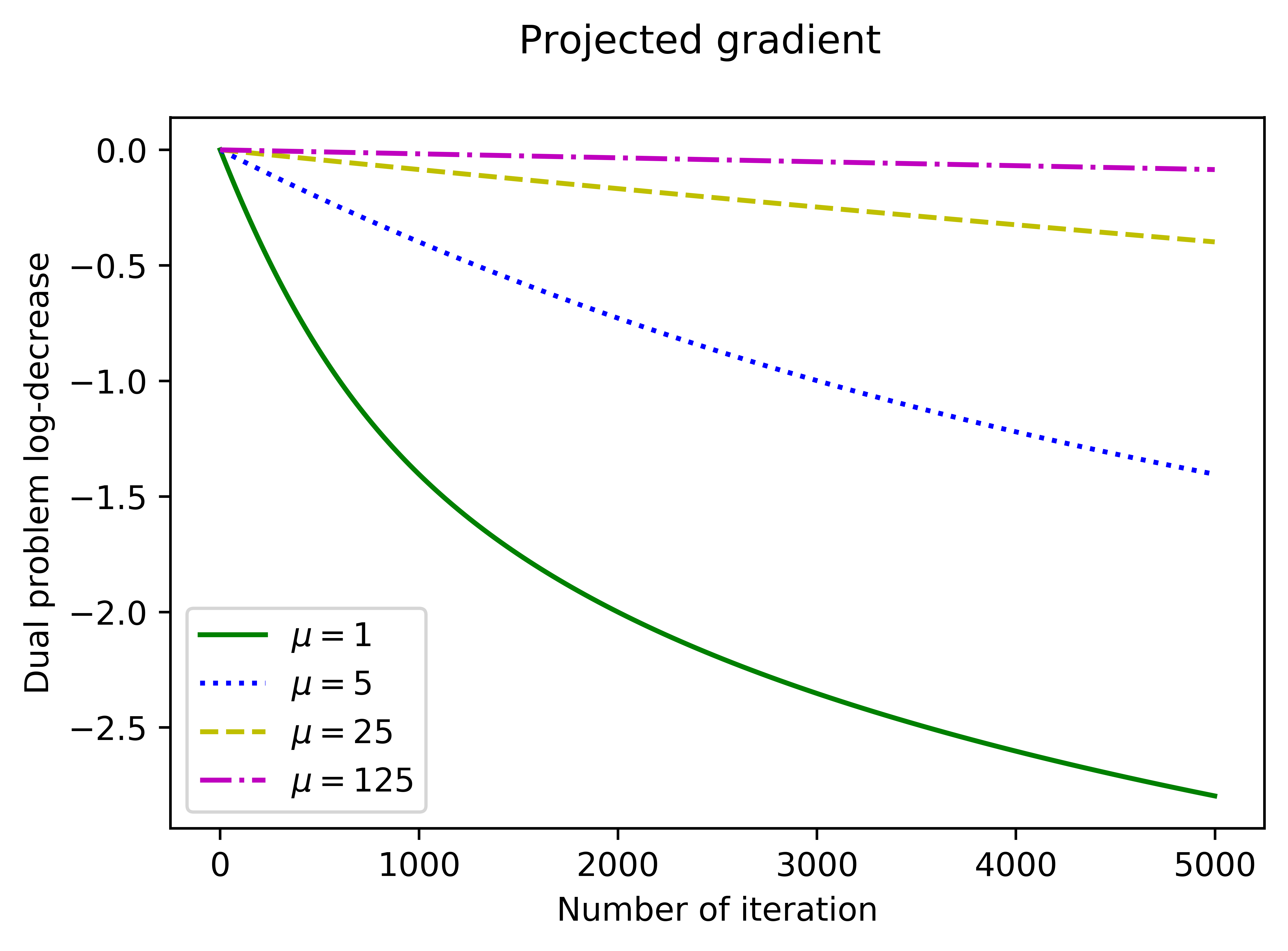}
    \includegraphics[width=0.4\textwidth]{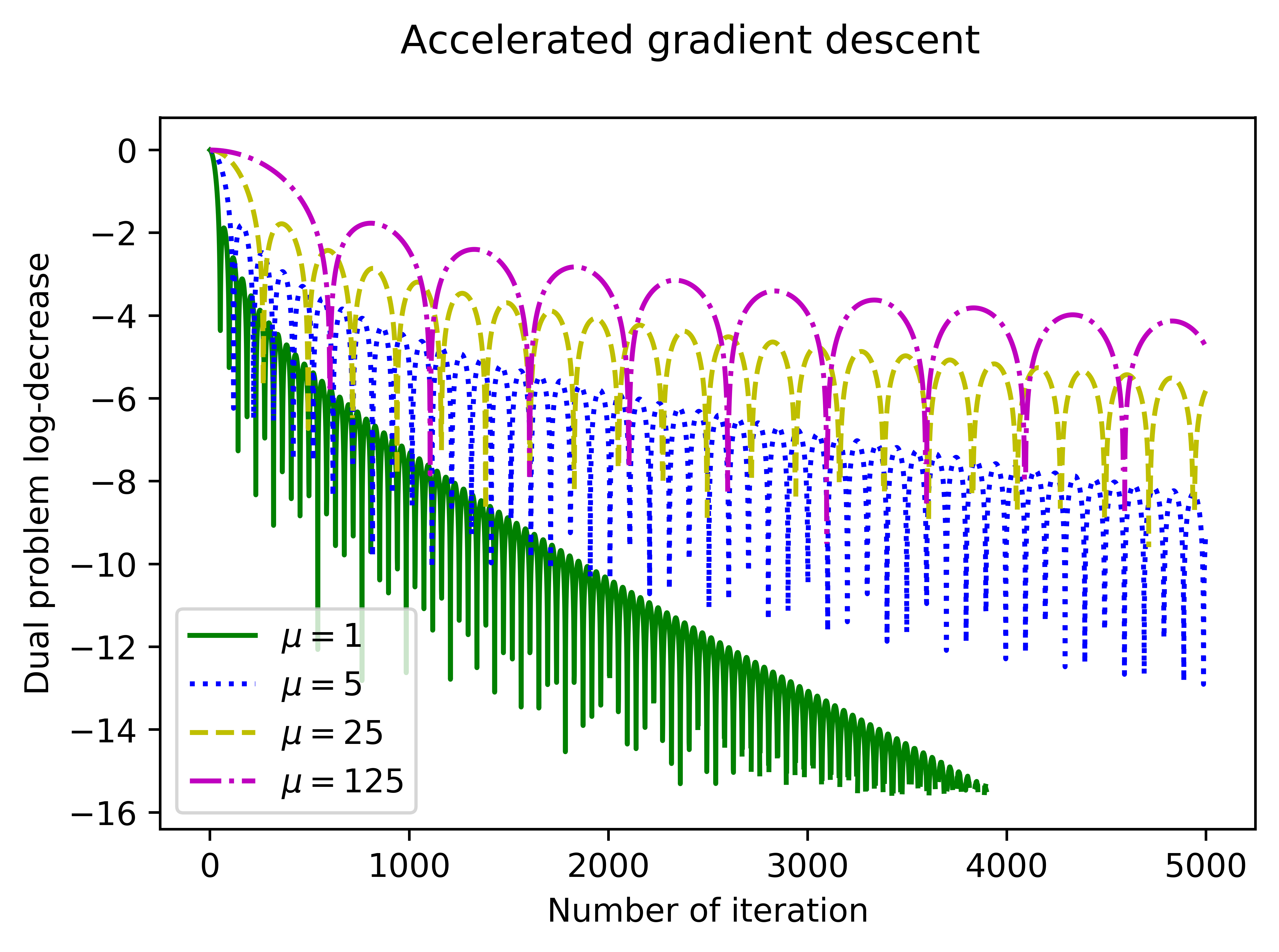}
\caption{Projected gradient vs Accelerated gradient descent,
 $100$~factories, cost functions~\eqref{eq_p5_prob1_fc1}, Dual value in log scale, mean of 20 run}
\label{mur_graph1}
\end{figure}

\begin{figure}[H]    \includegraphics[width=0.4\textwidth]{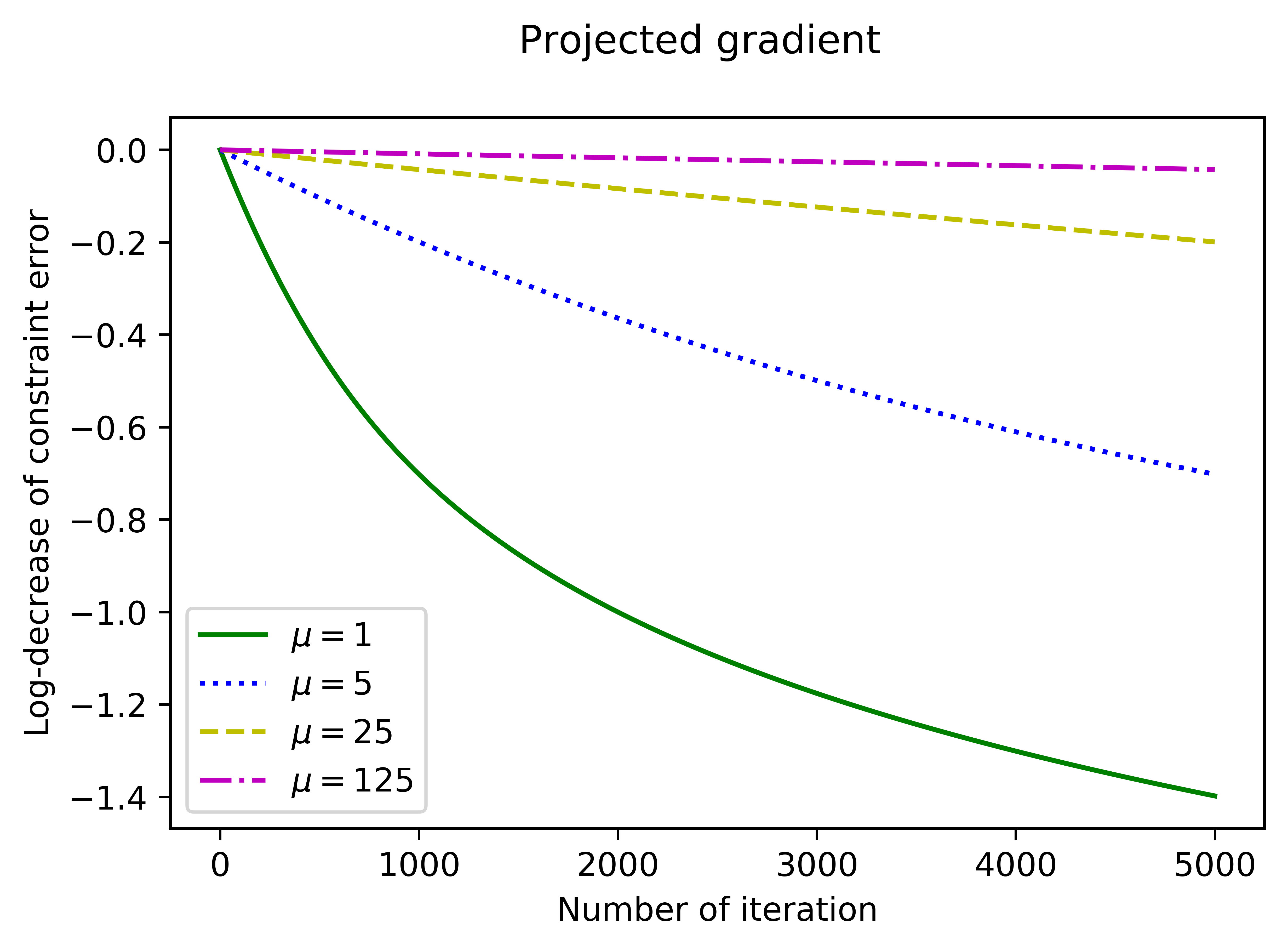}
    \includegraphics[width=0.4\textwidth]{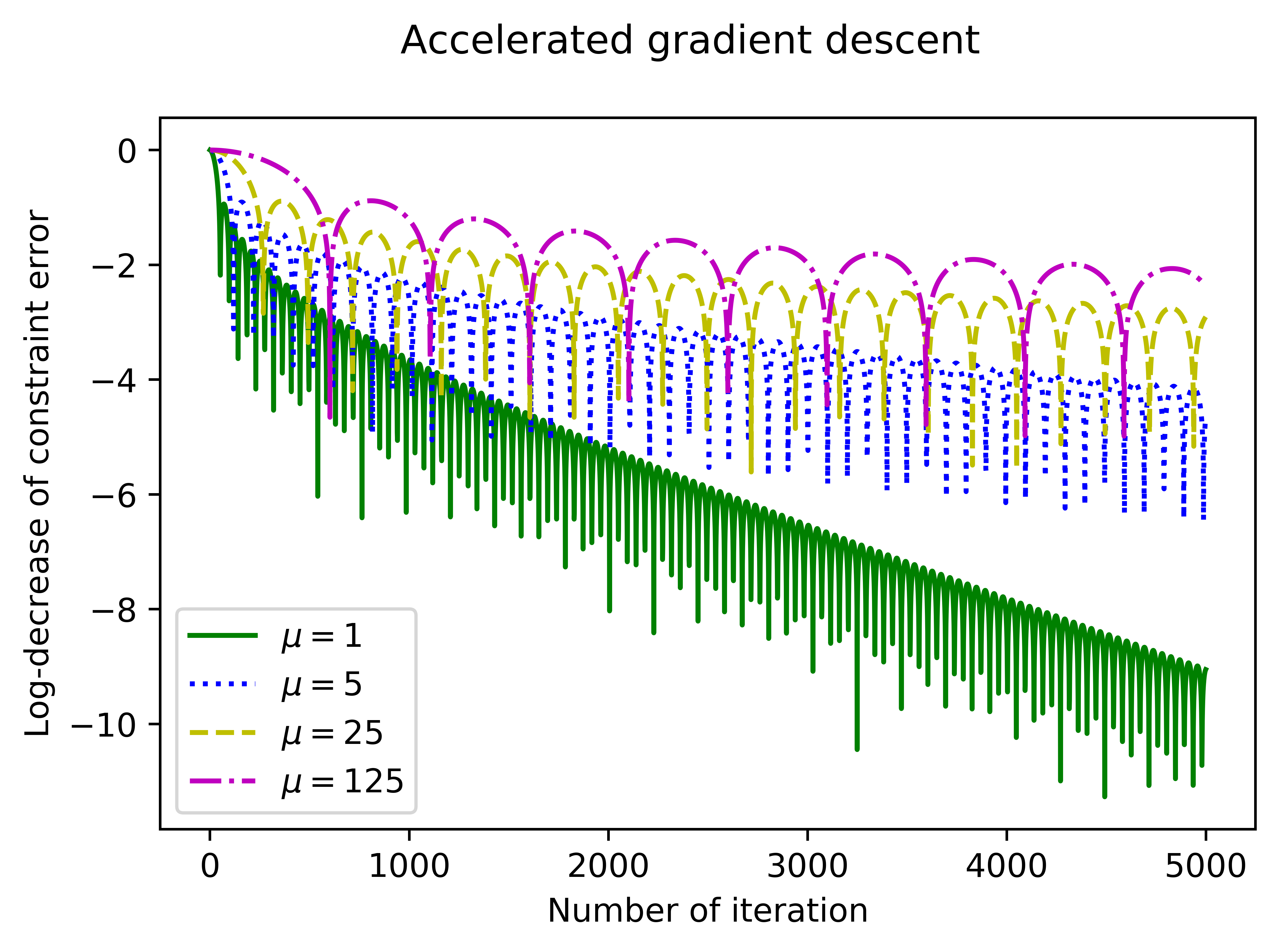}
\caption{Projected gradient vs Accelerated gradient descent,
 $100$~factories, cost functions~\eqref{eq_p5_prob1_fc1}, Inexactness of constraints in log scale, mean of 20 run}
\label{mur_graph2}
\end{figure}

\begin{figure}[H]    \includegraphics[width=0.4\textwidth]{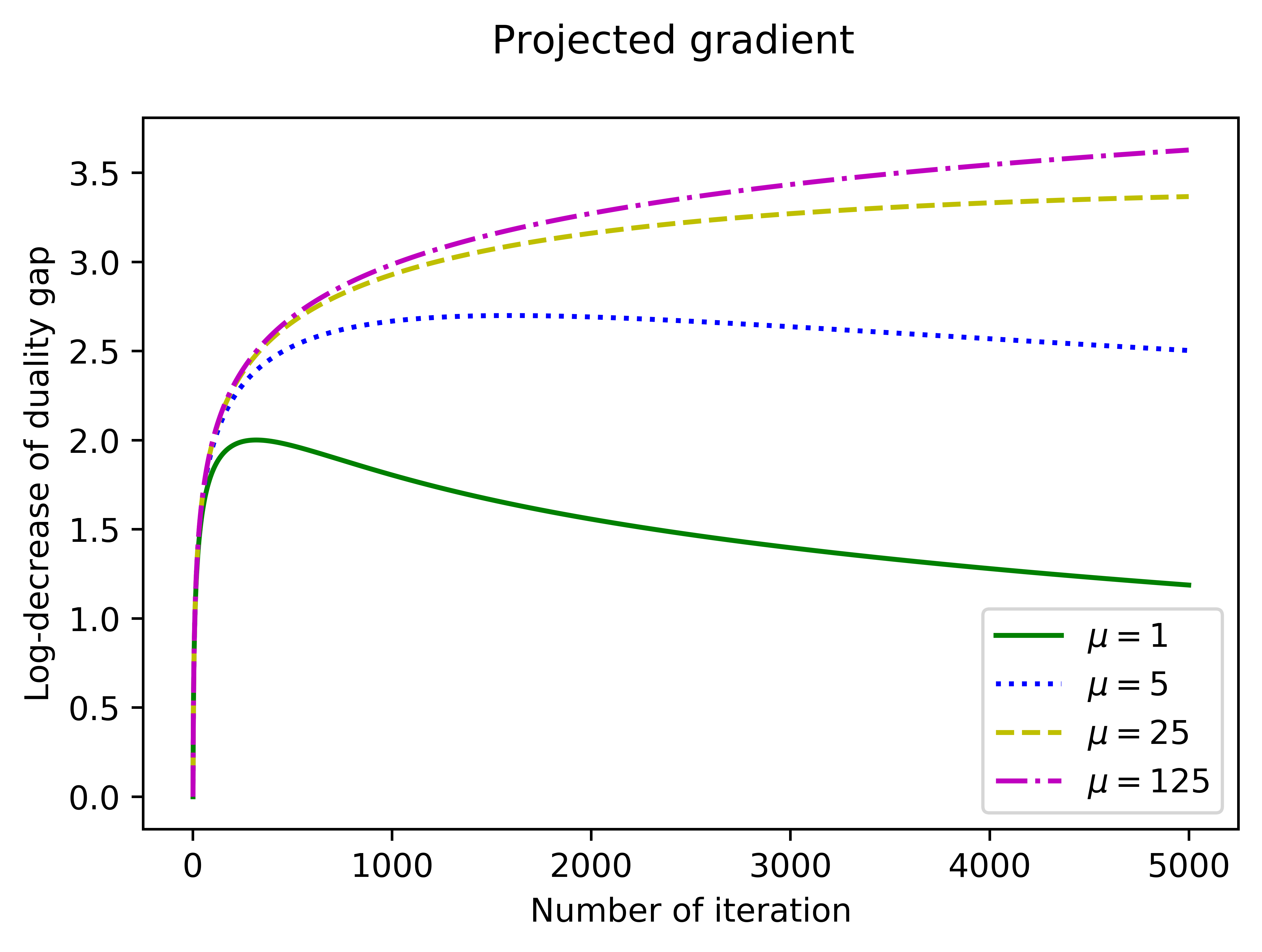}
    \includegraphics[width=0.4\textwidth]{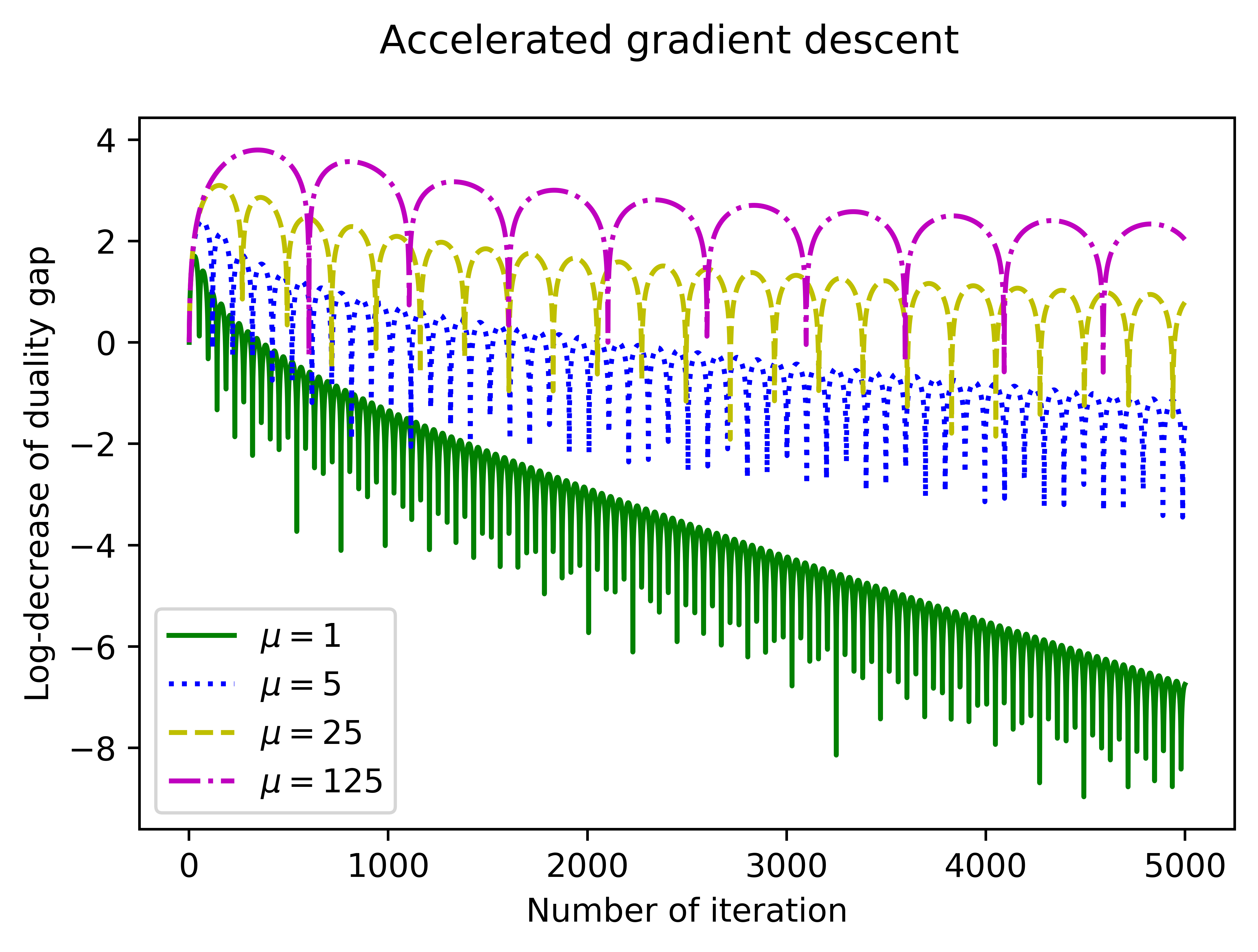}
\caption{Projected gradient vs Accelerated gradient descent,
 $100$~factories, cost functions~\eqref{eq_p5_prob1_fc1}, Duality gap value in log scale, mean of 20 run}
\label{mur_graph3}
\end{figure}

Secondly, we consider, how the results depend on the strong convexity parameter $\mu$. For that, we consider functions \eqref{eq_p5_prob1_fc1} with $\alpha \in N [100;400]$ and $\mu \in \{1, 5, 25, 125\}$. We make 20 runs with random $\alpha$ and show the mean of this runs.
Next in graphics the left side will be experiments for Projected gradient and the right side will be experiments for Accelerated gradient descent. The practical convergence rate for the objective value of the dual problem (D) is shown in log scale in fig.~\ref{mur_graph1}. The convergence rate for the constraint feasibility is shown in log scale in fig.~\ref{mur_graph2}.
The convergence rate for the value of the duality gap is shown in log scale  in fig.~\ref{mur_graph3}.

These plots show that the theoretical dependence on $\mu$ in both theorems is correct and occurs in practice.

\section*{Conclusion}

In this paper, we considered the resource allocation problem as a convex optimization problem. To solve this problem, we use and interpret gradient and accelerated gradient descent applied to the special variant of dual problem. We obtain
the convergence rate both for the primal iterates and the dual iterates. We obtain faster convergence rates than the ones known in the literature (for interpretable methods) and generalize our methods for the case of multiple producers and products.

The work of A.V. Gasnikov and P.E. Dvurechensky was supported by RFBR grant 18-29-03071 mk. The work of A.V. Gasnikov partially supported by RFBR 18-31-20005 mol-a-ved. The work of A.S. Ivanova was supported by Russian president grant MD-1320.2018.1. The work of D. Kamzolov was supported by Russian National Fund grant RNF 17-11-01027.

We'd like to thank Yu. Nesterov for the fruitful discussion.

\bibliographystyle{plain} 
\bibliography{gOMSguide}

\newpage

\section{Appendix}

\textbf{Proof of Lemma~
\ref{Lm:1}.}

Set $\bar x_1=\, \ldots \,=\bar x_n=\dfrac{2C}{n}$  and $\bar y_1=\, \ldots \,=\bar y_n=\dfrac{C}{n}$. And define $\bar {\mathtt{x}} = (\bar x_1, \, \ldots,\, \bar x_n)^{\top}$ and $\bar \y = (\bar y_1, \, \ldots,\, \bar y_n)^{\top}$. Notice that the point $\left(\bar {\mathtt{x}}^{\top}, \bar \y^{\top}\right)^{\top}$  satisfies the Slater's condition as 
$$
\dfrac{2C}{n}=\bar x_k > \bar y_k=\dfrac{C}{n}, \, \; \; k=1, \;  \ldots, \,n.
$$
Because the cost functions $f_k(x_k), \;k=1, \;  \ldots, \,n$ are increasing due to the economic interpretation, we obtain 

\begin{equation*}
\begin{split}
\sum \limits_{k=1}^n f_k(0) & =  \min \limits_{\substack{\sum \limits_{k=1}^n y_k \geqslant  C, \;  y_k \geqslant  0,\;\\ x_k \geqslant  0,\; k=1, \, \ldots, \, n}}\sum \limits_{k=1}^n f_k(x_k)\\
&=\min \limits_{\substack{\sum \limits_{k=1}^n y_k \geqslant  C, \;  y_k \geqslant  0,\;\\ x_k \geqslant  0,\; k=1, \, \ldots, \, n}}\left\{\sum \limits_{k=1}^n f_k(x_k)+\sum \limits_{k=1}^n\left(y_k-x_k\right)\underbrace{p_k}_{=\, 0}\right\}\\
&\leqslant   \max \limits_{p \geqslant  0}\min \limits_{\substack{\sum \limits_{k=1}^n y_k \geqslant  C, \;  y_k \geqslant  0,\;\\ x_k \geqslant  0,\; k=1, \, \ldots, \, n}}\left\{\sum \limits_{k=1}^n f_k(x_k)+\sum \limits_{k=1}^n\left(y_k-x_k\right)p_k\right\}\\
&= \min \limits_{\substack{\sum \limits_{k=1}^n y_k \geqslant  C, \;  y_k \geqslant  0,\;\\ x_k \geqslant  0,\; k=1, \, \ldots, \, n}}\left\{\sum \limits_{k=1}^n f_k(x_k)+\sum \limits_{k=1}^n\left(y_k-x_k\right)p_k^{*}\right\}\\
&\leqslant  \sum \limits_{k=1}^n f_k(\bar x_k)+\sum \limits_{k=1}^n\left(\bar y_k-\bar x_k\right)p^{*}_k\\
&\leqslant    \sum \limits_{k=1}^n f_k(\bar x_k)-\dfrac{C}{n}\sum \limits_{k=1}^{n}p^{*}_k.
\end{split}
\end{equation*}
Since $p_k^{*}\geqslant 0,\; k=1, \, \ldots, \, n$ we obtain that  
\begin{equation*}
    \norm{\p^{*}}_{1}\leqslant p_{max},
\end{equation*}
where 
\begin{equation*}
   p_{max}= \dfrac{n}{C}\left(\sum \limits_{k=1}^n f_k\left(\dfrac{2C}{n}\right)-\sum \limits_{k=1}^n f_k(0)\right). 
\end{equation*}
And since that we obtain the upper bound for each component 
$$
p^*_k\leqslant p_{max}, \, \, k=1, \, \ldots, \, n
$$
from which we obtain the statement of the Lemma.
\qed

\textbf{Proof of Lemma~
\ref{lemma2}.}

Using~\eqref{tildephi}, we can rewrite the step as 
\begin{equation*}
\begin{split}
    \p^{t+1}&=\argmin_{\p\geqslant 0}\left\{C\max \limits_{k=1, \, \ldots, \, n}(-p_k)+\left<\mathtt{x}(\p^t),\p-\p^{t}\right>+\dfrac{L_{\psi}}{2}\Big(\norm{\p}_2^2-2\left<\p,\p^{t}\right> +\norm{\p^t}_{2}^{2}\Big)\right\}\\&
    =\argmin_{\p\geqslant 0}\left\{C\max \limits_{k=1, \, \ldots, \, n}(-p_k)+\dfrac{L_{\psi}}{2}\Big(\norm{\p}_2^2-2\left<\p,\p^{t}\right>+\dfrac{2}{L_{\psi}}\left<\mathtt{x}(\p^t),\p-\p^{t}\right> +\norm{\p^t}_{2}^{2}\Big)\right\}\\&
    =\argmin_{\p\geqslant 0}\left\{\max \limits_{k=1, \, \ldots, \, n}(-p_k)+\dfrac{L_{\psi}}{2C}\norm{\p-\left(\p^t-\dfrac{1}{L_{\psi}}\mathtt{x}(\p^t) \right)}_{2}^{2}\right\}.
\end{split}
\end{equation*}
Let us define $\gamma=\dfrac{C}{L_{\psi}}$. Then to determine $\p^{t+1}$ it is necessary to solve the following problem
$$
\max\limits_{k=1,\, \ldots\, n}(-p_k)+\dfrac{1}{2\gamma}\norm{\p-\tilde{\p}^{t+1}}^2_2 \rightarrow \min\limits_{p_k\geqslant 0,\, k=1,\, \ldots\, n,} 
$$
which we can rewrite as the following equivalent problem
$$
\eta_{t+1}+\dfrac{1}{2\gamma}\norm{\p-\tilde{\p}^{t+1}}^2_2 \rightarrow \min\limits_{ \eta_{t+1} \leqslant 0,\, \eta_{t+1} \geqslant -p_k, \; k=1,\, \ldots \,n}
$$
where $\eta_{t+1} \in \R$. The Lagrangian is 
$$
L(\p,\eta_{t+1},\z,w)=\eta_{t+1}+\dfrac{1}{2\gamma}\norm{\p-\tilde{\p}^{t+1}}^2_2+\sum_{k=1}^n z_k(-p_k-\eta_{t+1})+w\eta_{t+1}
$$
with dual variable $\z \in \R^{n}_{+}$ and $w \in \R_{+}$. Let $(\p^{*},\eta^{*}_{t+1},\z^{*},w^{*})$ be a solution, then the optimality conditions are:
\begin{align}
    \label{eq1}
    -p^{*}_k\leqslant \eta^{*}_{t+1}&, \, k=1, \, \ldots \,,n.\\
    \label{eq2}
    z^{*}_k\geqslant 0&, \, k=1, \, \ldots \,,n.\\
\label{eq3}
   z_k^{*}(p^{*}_k+\eta^{*}_{t+1})=0&, \, k=1, \, \ldots \,,n.\\
\label{eq4}
    \dfrac{1}{\gamma}(p^{*}_k-\tilde{p}_k^{t+1})-z^{*}_k=0&, \, k=1, \, \ldots \,,n.\\
\label{eq5}
    \sum_{k=1}^n z^*_k&=1+w^{*}.\\
\label{eq6}
    w^{*}\eta_{t+1}^{*}&=0.\\
\label{eq7}
     w^{*}&\geqslant 0.\\
\label{eq8}
     \eta^*_{t+1}&\leqslant 0.
\end{align}
\begin{itemize}
    \item If $\eta^{*}_{t+1}=0$.
    \begin{itemize}
        \item If $\eta_{t+1}^{*}=-p_k^{*}=0$, then from~\eqref{eq4} and \eqref{eq2} obtain, that 
        \begin{equation}
\label{eq4.1}
     z_{k}^{*}=\dfrac{1}{\gamma}(-\tilde{p}^{t+1}_k)=\dfrac{1}{\gamma}(-\tilde{p}^{t+1}_k)_{+}.
\end{equation}
        \item If $p_{k}^{*}>-\eta^*_{t+1}=0$ then from~\eqref{eq3} and ~\eqref{eq4}
        $$ 
        z^{*}_{k}=0,\, \,  p_k^*=\tilde{p}^{t+1}_k.
        $$
Also since in this case $p^{*}_k>0$, then $\tilde{p}^{t+1}_k>0$. And using this we obtain
\begin{equation}
\label{eq4.2}
     z_{k}^{*}=0=\dfrac{1}{\gamma}(-\tilde{p}^{t+1}_k)_{+}.
\end{equation}
    \end{itemize}
Then, using~\eqref{eq4.1}, ~\eqref{eq4.2}, ~\eqref{eq5}  and \eqref{eq7} obtain, that 
$$
\sum_{k=1}^{n}(-\tilde{p}^{t+1}_k)_{+}=\gamma(1+w^*)\geqslant \gamma.
$$
As a result, we obtain that for $\sum\limits_{k=1}^{n}(-\tilde{p}^{t+1}_k)_{+} \geqslant \gamma$ the case then $\eta^{*}_{t+1}=0$ and solution determines as 
$$
p_k^{*}=\max\{0,\tilde{p}^{t+1}_k\}.
$$
\item If $\eta_{t+1}^{*}<0$ then $w^*=0$.
\begin{itemize}
    \item If $p_k^{*}=-\eta_{t+1}^{*}$ then from~\eqref{eq4} obtain
    \begin{equation}
\label{eq4.3}
     z_{k}^{*}=\dfrac{1}{\gamma}(-\eta_{t+1}^{*}-\tilde{p}^{t+1}_k)_{+}.
\end{equation}
    \item If $p_k^{*}>-\eta_{t+1}^{*}$ then from~\eqref{eq3} and ~\eqref{eq4} we obtain that $z^{*}_{k}=0$, $p_{k}^{*}=\tilde{p}^{t+1}_{k}$. Also since in this case $\tilde{p}^{t+1}_k=p_{k}^{*}>-\eta_{t+1}^{*}$, then 
    \begin{equation}
    \label{eq4.4}
     z_{k}^{*}=0=\dfrac{1}{\gamma}(-\eta_{t+1}^{*}-\tilde{p}^{t+1}_k)_{+}.
\end{equation}
\end{itemize}
Finally using~\eqref{eq4.3}, ~\eqref{eq4.4} and~\eqref{eq5} obtain that
$$
\sum_{k=1}^{n}(-\eta_{t+1}^{*}-\tilde{p}^{t+1}_k)_{+}=\gamma.
$$
and solution determines as  
$$
p_k^{*}=\max\{-\eta_{t+1}^{*},\,\tilde{p}^{t+1}_k\}.
$$
\end{itemize} 
So, define $p^{t+1}_{center} =-\eta_{t+1}^{*}$ we finishes the proof.
\qed

\textbf{Proof of Theorem~
\ref{th2}.}

Let's consider the step 
\begin{equation*}
    \p^{t+1}=\argmin_{\p\geqslant 0}\underbrace{\left\{\left<\nabla\psi(\p^t),\p-\p^{t}\right>-C\min \limits_{k=1, \, \ldots, \, n}p_k+\dfrac{L_{\psi}}{2}\norm{\p-\p^t}_{2}^{2}\right\}}_{\Phi(\p,\p^t)}.
\end{equation*}
The function $\Phi(\p,\p^t)$ is convex with respect to the variable $\p$, thеn from the definition of $\p^{t+1}$ follows that 
\begin{equation}
\label{min_cond}
    \left<\nabla_{\p^{t+1}}\Phi(\p^{t+1},\p^t),\p-\p^{t+1}\right> \geqslant 0, \, \, \forall \p \geqslant 0.
\end{equation}
Define function 
$$
\tilde{\Phi}(\p,\p^t)=\left<\nabla\psi(\p^t),\p-\p^{t}\right>-C\min \limits_{k=1, \, \ldots, \, n}p_k =\left<\nabla\psi(\p^t),\p-\p^{t}\right> + g(\p),
$$
which is convex with respect to the variable $\p$. From~\eqref{min_cond} we obtain
\begin{align}
 0& \leqslant \left<\nabla_{\p^{t+1}}\Phi(\p^{t+1},\p^t),\p-\p^{t+1}\right> = \left<\nabla_{\p^{t+1}}\tilde{\Phi}(\p^{t+1},\p^t)+L_{\psi}(\p^{t+1}-\p^{t}),\p-\p^{t+1}\right> \notag \\& \overset{conv-ty}{\leqslant} \tilde{\Phi}(\p,\p^t)-\tilde{\Phi}(\p^{t+1},\p^t)+L_{\psi}\left<\p^{t+1}-\p^{t},\p-\p^{t+1}\right>.
 \label{eq2.6}
\end{align}
Notice, that 
\begin{align}
L_{\psi}\left<\p^{t+1}-\p^{t},\p-\p^{t+1}\right>&=L_{\psi}\left<\p^{t+1},\p\right>+L_{\psi}\left<\p^{t},\p^{t+1}\right>-L_{\psi}\left<\p^{t+1},\p^{t+1}\right>-L_{\psi}\left<\p^{t},\p\right> \notag \\&
=\dfrac{L}{2}\left<\p,\p\right>-L_{\psi}\left<\p^{t},\p\right>+\dfrac{L_{\psi}}{2}\left<\p^{t},\p^{t}\right>-\dfrac{L_{\psi}}{2}\left<\p,\p\right>+L_{\psi}\left<\p,\p^{t+1}\right> \notag \\&
-\dfrac{L_{\psi}}{2}\left<\p^{t+1},\p^{t+1}\right>-\dfrac{L_{\psi}}{2}\left<\p^{t+1},\p^{t+1}\right>-\dfrac{L_{\psi}}{2}\left<\p^{t},\p^{t}\right>+L_{\psi}\left<\p^{t},\p^{t+1}\right> \notag \\&
=\dfrac{L_{\psi}}{2}\norm{\p-\p^t}_2^2-\dfrac{L_{\psi}}{2}\norm{\p-\p^{t+1}}_2^2-\dfrac{L_{\psi}}{2}\norm{\p^{t+1}-\p^{t}}_2^2.
\label{div_eq}
\end{align}
Substituting this inequality to~\eqref{eq2.6}, we obtain 
\begin{equation*}
\begin{split}
0&\leqslant \tilde{\Phi}(\p,\p^t)-\tilde{\Phi}(\p^{t+1},\p^t)+\dfrac{L_{\psi}}{2}\norm{\p-\p^t}_2^2-\dfrac{L_{\psi}}{2}\norm{\p-\p^{t+1}}_2^2-\dfrac{L_{\psi}}{2}\norm{\p^{t+1}-\p^{t}}_2^2\\ &
=\left<\nabla\psi(\p^t),\p-\p^{t}\right> + g(\p)-\left<\nabla\psi(\p^t),\p^{t+1}-\p^{t}\right> - g(\p^{t+1})+\dfrac{L_{\psi}}{2}\norm{\p-\p^t}_2^2\\&
-\dfrac{L_{\psi}}{2}\norm{\p-\p^{t+1}}_2^2-\dfrac{L_{\psi}}{2}\norm{\p^{t+1}-\p^{t}}_2^2
\end{split}
\end{equation*}
Finally, since~\eqref{Lips_grad} we can rewrite as 
\begin{equation}
\label{Lip_cond}
\psi(\p^{1}) \leqslant \psi(\p^2)+ \left<\nabla\psi(\p^2),\p^{1}-\p^{2}\right>+\dfrac{L_{\psi}}{2}\norm{\p^{1}-\p^{2}}_2^2, \, \forall \p^1, \, \p^2 \geqslant \mathtt{0}
\end{equation}
we obtain that
$$
\psi(\p^{t+1}) \leqslant \psi(\p^t)+\left<\nabla\psi(\p^t),\p-\p^{t}\right> + g(\p)- g(\p^{t+1})+\dfrac{L_{\psi}}{2}\norm{\p-\p^t}_2^2
-\dfrac{L_{\psi}}{2}\norm{\p-\p^{t+1}}_2^2
$$
or 
$$
\varphi(\p^{t+1}) \leqslant \varphi(\p^t)+\left<\nabla\psi(\p^t),\p-\p^{t}\right>-g(\p^{t}) + g(\p)+\dfrac{L_{\psi}}{2}\norm{\p-\p^t}_2^2
-\dfrac{L_{\psi}}{2}\norm{\p-\p^{t+1}}_2^2.
$$
Summing these inequalities from $t=0$ to $t=N-1$ we get, using
the convexity of $\varphi$ and that it is true for all $\p \geqslant \mathtt{0}$ 
\begin{equation}
\label{fin_est1}
    \varphi(\p^{N})\leqslant\dfrac{1}{N}\min\limits_{p \geqslant \mathtt{0}}\left\{\sum\limits_{t=0}^{N-1}\left[\varphi(\p^t)+\left<\nabla\psi(\p^t),\p-\p^t\right>-g(\p^t)\right]+g(\p)+\dfrac{L_{\psi}}{2}\norm{\p-\p^0}_{2}^{2}\right\},
\end{equation}
where $\p^{N}=\dfrac{1}{N}\sum\limits_{t=1}^{N}\p^{t}$. Let us define the starting point $\p^0$ satisfying $0 \leqslant p_k^0 \leqslant p_{max},\; k=1, \, \ldots, \, n$, where $p_{max}$ is given in \eqref{eq:p_max}. Then we obtain that
\begin{equation*}
   \norm{\p^{0}}_{2}\leqslant \sqrt{n}p_{max}.
\end{equation*}
Let us introduce a set
$$
B^{+}_{3R}(0)=\{\mathtt{p}\;:\; \mathtt{p}\geqslant  0, \; \left|\left|\mathtt{p}-\mathtt{p}^0\right|\right|_2\leqslant   3R\},
$$
where
\begin{equation}
\label{4_R}
\left|\left|\mathtt{p}^0-\mathtt{p}^{*}\right|\right|_2+\norm{\p^0}_2\leqslant 2\norm{\p^0}_2+\norm{\p^*}_2=3p_{max}\sqrt{n}=R,
\end{equation}
herewith all the obtaining $\mathtt{p}^t$ will consist in $B^{+}_{2R}(\mathtt{0})$ :
\begin{equation}
\label{bound_2R}
\|\mathtt{p}^t\|_2 \leqslant 2R,
\end{equation}
since (second paragraph \cite{Gasnikov})
\begin{align*}
\|\mathtt{p}^t\|_2 = \|\mathtt{p}^t - \mathtt{p}^0\|_2 +\norm{\p^0}_2&\leqslant \|\mathtt{p}^t - \mathtt{p}^{*}\|_2 + \|\mathtt{p}^{*} - \mathtt{p}^0\|_2+\norm{\p^0}_2 \\& \leqslant 2 \|\mathtt{p}^{*} - \mathtt{p}^0\|_2+\norm{\p^0}_2 = 2 \|\mathtt{p}^{*}\|_2+3\norm{\p^0}_2 = 5p_{max}\sqrt{n}\leqslant 2R.
\end{align*}
Since that we can rewrite~\eqref{fin_est1} as 
\begin{equation*}
\begin{split}
    \varphi(\p^{N})&\leqslant\dfrac{1}{N}\min\limits_{\p \geqslant \mathtt{0}}\left\{\sum\limits_{t=0}^{N-1}\left[\varphi(\p^t)+\left<\nabla\psi(\p^t),\p-\p^t\right>-g(\p^t)\right]+g(\p)+\dfrac{L_{\psi}}{2}\norm{\p-\p^0}_{2}^{2}\right\}\\& \leqslant
    \dfrac{1}{N}\min\limits_{\p\in B_{3R}^{+}(\mathtt{0})}\left\{\sum\limits_{t=0}^{N-1}\left[\psi(\p^t)+\left<\nabla\psi(\p^t),\p-\p^t\right>\right]+g(\p)+\dfrac{L_{\psi}}{2}\norm{\p-\p^0}_{2}^{2}\right\}.
\end{split}
\end{equation*}

Notice that the minimum is taken over the set $B_{3R}^{+}(\mathtt{0})=B_{3R}^{+}(\p^0)$ and not over the set $B_{3R}^{+}(\p^{*})$, since the $\p^*$ is not known and, according to the definition of $R$, $\p^* \in B_{R}^{+}(\p^0).$ 
Define 
\begin{equation*}
    \tilde{\p}=\argmin\limits_{\p \in B_{3R}^{+}(\mathtt{0})}\left\{\sum\limits_{t=0}^{N-1}\left[\psi(\p^t)+\left<\nabla\psi(\p^t),\p-\p^t\right>\right]+g(\p) \right\},
\end{equation*}
then, since
$$
\norm{\tilde{\p}-\p^0}_{2}^{2} \leqslant 2\norm{\tilde{\p}}_{2}^{2}+2\norm{\p^0}_{2}^{2}\leqslant 18R^2+\dfrac{2}{9}R^2=\dfrac{164}{9}R^2
$$
we obtain 
\begin{equation*}
\begin{split}
    \varphi(\p^{N})&\leqslant
    \dfrac{1}{N}\min\limits_{\p \in B_{3R}^{+}(\mathtt{0})}\left\{\sum\limits_{t=0}^{N-1}\left[\psi(\p^t)+\left<\nabla\psi(\p^t),\p-\p^t\right>\right]+g(\p)+\dfrac{L_{\psi}}{2}\norm{\p-\p^0}_{2}^{2}\right\}\\&
    \leqslant
    \dfrac{1}{N}\min\limits_{\p \in B_{3R}^{+}(\mathtt{0})}\left\{\sum\limits_{t=0}^{N-1}\left[\psi(\p^t)+\left<\nabla\psi(\p^t),\p-\p^t\right>\right]+g(\p)\right\}+\dfrac{L_{\psi}}{2N}\norm{\tilde{\p}-\p^0}_{2}^{2}
    \\&
    \leqslant
    \dfrac{1}{N}\min\limits_{\p \in B_{3R}^{+}(\mathtt{0})}\left\{\sum\limits_{t=0}^{N-1}\left[\psi(\p^t)+\left<\nabla\psi(\p^t),\p-\p^t\right>\right]+g(\p)\right\}+\dfrac{82L_{\psi}R^2}{9N}.
\end{split}
\end{equation*}
Notice, that 
\begin{equation*}
\begin{split}
   &-\dfrac{1}{N}\min\limits_{\p \in B_{3R}^{+}(\mathtt{0})}\left\{\sum\limits_{t=0}^{N-1}\left[\psi(\p^t)+\left<\nabla\psi(\p^t),\p-\p^t\right>\right]+g(\p)\right\} \\& =\dfrac{1}{N}\max\limits_{\p \in B_{3R}^{+}(\mathtt{0})}\left\{\sum\limits_{t=0}^{N-1}\left[ \left<\mathtt{x}(\p^t),\p^t-\p\right> + \sum\limits_{k=1}^{n} f_{k}(\mathtt{x}_{k}(\p^{t}))- \left<\mathtt{x}(\p^t),\p^t\right>\right] - g(\p) \right\}
   \\&= \dfrac{1}{N} \sum\limits_{t=0}^{N-1} f(\mathtt{x}(\p^{t})) + \dfrac{1}{N}\max\limits_{\p \in B_{3R}^{+}(\mathtt{0})}\left\{\sum\limits_{t=0}^{N-1}\left<- \mathtt{x}(\p^t),\p\right> - g(\p) \right\} \\&\overset{conv-ty}{\geqslant} f(\mathtt{x}^{N}) + \dfrac{1}{N}\max\limits_{\p \in B_{3R}^{+}(\mathtt{0})}\left\{\sum\limits_{t=0}^{N-1} \left<- \mathtt{x}(\p^t),\p\right> - g(\p) \right\}
\\&= f(\mathtt{x}^{N}) + \max\limits_{\p \in B_{3R}^{+}(\mathtt{0})}\left\{ - \left< \mathtt{x}^N,\p\right>+C\min\limits_{k=1, \, \dots \, ,n}p_{k} \right\},
\end{split}
\end{equation*}
where $\mathtt{x}^{N}=\dfrac{1}{N}\sum\limits_{t=0}^{N-1}\mathtt{x}(\p^{t})$. Since 
$$ 
C\min\limits_{k=1, \, \ldots \, ,n}p_{k}-  \sum\limits_{k=1}^{n}x_k^N \min\limits_{k=1, \, \ldots \, ,n}p_{k} \geqslant C\min\limits_{k=1, \, \ldots \, ,n}p_{k}-  \sum\limits_{k=1}^{n}x_k^N p_{k},
$$
we obtain 
\begin{equation*}
\begin{split}
 \max\limits_{\p \in B_{3R}^{+}(\mathtt{0})}\left\{ - \left< \mathtt{x}^N,\p\right>  +C\min\limits_{k=1, \, \ldots \, ,n}p_{k} \right\}
&= \max\limits_{\substack{\p \in B_{3R}^{+}(\mathtt{0})\;\\
\\p_1=\ldots=p_n=p}}\left\{ - \left< \mathtt{x}^N,\p\right>+C\min\limits_{k=1, \, \ldots \, ,n}p_{k} \right\}\\&
=
\max\limits_{\substack{\p \in B_{3R}^{+}(\mathtt{0})\;\\
\\p_1=\ldots=p_n=p}}\left\{Cp- \sum\limits_{k=1}^{n}x_k^N p \right\}
=\dfrac{3R}{\sqrt{n}}\left[C- \sum\limits_{k=1}^{n}x^N_k \right]_{+}.
\end{split}
\end{equation*}
Thus, we obtain the following estimation 
\begin{equation}
\label{estimation}
    \varphi(\p^{N})+f(\mathtt{x}^{N})+\dfrac{3R}{\sqrt{n}}\left[C- \sum\limits_{k=1}^{n}x_k^N \right]_{+}\leqslant \dfrac{82L_{\psi}R^2}{9N}. 
\end{equation}
Considering the weak duality $-f(\mathtt{x}^*)\leqslant\varphi(\p^*)$, obtain 
\begin{equation*}
\begin{split}
f(\mathtt{x}^N)-f(\mathtt{x}^{*})&\leqslant   f(\mathtt{x}^N)+\varphi(\mathtt{p}^{*})\leqslant   \varphi(\mathtt{p}^N)+f(\mathtt{x}^N)\\
&\leqslant   \varphi(\mathtt{p}^N)+f(\mathtt{x}^N)+\dfrac{3R}{\sqrt{n}}\left[C- \sum\limits_{k=1}^{n}x_k^N \right]_{+}\leqslant \dfrac{82L_{\psi}R^2}{9N}=\varepsilon.
\end{split}   
\end{equation*}
Using~\eqref{bound_2R},~\eqref{estimation} and $(D)$ we get 
\begin{align}
\dfrac{R}{\sqrt{n}}\left[C- \sum\limits_{k=1}^{n}x_k^N \right]_{+}&\leqslant   \overbrace{\underbrace{\left<\mathtt{x}^{N}, \, \mathtt{p}^N\right> - C \min\limits_{k=1, \, \ldots \,,n}p^{N}_{k} }_{\geqslant  - \dfrac{2R}{\sqrt{n}}\left[C- \sum\limits_{k=1}^{n}x_k^N \right]_{+}}-f(\mathtt{x}^N)}^{\leqslant   \varphi(\mathtt{p}^N)}+f(\mathtt{x}^N) 
+\dfrac{3R}{\sqrt{n}}\left[C- \sum\limits_{k=1}^{n}x_k^N \right]_{+}\leqslant   \varepsilon. 
\label{duality_gap1}
\end{align}
And since $R=3p_{max}\sqrt{n}$ and $L_{\psi}=\dfrac{n}{\mu}$ we obtain the statement of the theorem.

\qed

\textbf{Proof of Theorem~
\ref{th3}.}

Let's consider the step 
\begin{equation*}
    \mathtt{\y}^{t+1}=\argmin_{\p\geqslant 0}\underbrace{\left\{\alpha_{t+1}\left(\left<\nabla\psi(\p^{t+1}),\p-\p^{t+1}\right>-C\min \limits_{k=1, \, \ldots, \, n}p_k\right)+\dfrac{1}{2}\norm{\p-\y^t}_{2}^{2}\right\}}_{\Psi(\p,\p^{t+1})}.
\end{equation*}
The function $\Psi(\p,\p^t)$ is convex with respect to the variable $\p$, then from the definition of $\y^{t+1}$ follows that 
\begin{equation}
\label{min_cond_1}
    \left<\nabla_{\y^{t+1}}\Psi(\y^{t+1},\p^{t+1}),\p-\y^{t+1}\right> \geqslant 0, \, \, \forall \p \geqslant 0.
\end{equation}
Define function 
$$
\tilde{\Psi}(\p,\p^t)=\left<\nabla\psi(\p^t),\p-\p^{t}\right>-C\min \limits_{k=1, \, \ldots, \, n}p_k =\left<\nabla\psi(\p^t),\p-\p^{t}\right> + g(\p),
$$
which is convex with respect to the variable $\p$. From~\eqref{min_cond_1} we obtain
\begin{align}
 0& \leqslant \left<\nabla_{\y^{t+1}}\Psi(\y^{t+1},\p^{t+1}),\p-\y^{t+1}\right> = \left<\alpha_{t+1}\nabla_{\y^{t+1}}\tilde{\Psi}(\y^{t+1},\p^{t+1})+(\y^{t+1}-\y^{t}),\p-\y^{t+1}\right> \notag \\& \overset{conv-ty}{\leqslant} \alpha_{t+1}\tilde{\Psi}(\p,\p^{t+1})-\alpha_{t+1}\tilde{\Psi}(\y^{t+1},\p^{t+1})+\left<\y^{t+1}-\y^{t},\p-\y^{t+1}\right>.
 \label{eq3.6}
\end{align}
Using~\eqref{div_eq}, we obtain 
\begin{equation*}
\begin{split}
0&\leqslant \alpha_{t+1}\tilde{\Psi}(\p,\p^{t+1})-\alpha_{t+1}\tilde{\Phi}(\y^{t+1},\p^{t+1})+\dfrac{1}{2}\norm{\p-\y^t}_2^2-\dfrac{1}{2}\norm{\p-\y^{t+1}}_2^2-\dfrac{1}{2}\norm{\y^{t+1}-\y^{t}}_2^2\\ &
=\alpha_{t+1}\left<\nabla\psi(\p^{t+1}),\p-\p^{t+1}\right> + \alpha_{t+1}g(\p)-\alpha_{t+1}\left(\left<\nabla\psi(\p^{t+1}),\y^{t+1}-\p^{t+1}\right> + g(\y^{t+1})\right)\\&
+\dfrac{1}{2}\norm{\p-\y^t}_2^2
-\dfrac{1}{2}\norm{\p-\y^{t+1}}_2^2-\dfrac{1}{2}\norm{\y^{t+1}-\y^{t}}_2^2
\end{split}
\end{equation*}
Notice, that 
\begin{equation*}
\begin{split}
&-\alpha_{t+1}\left(\left<\nabla\psi(\p^{t+1}),\y^{t+1}-\p^{t+1}\right> + g(\y^{t+1})\right)-\dfrac{1}{2}\norm{\y^{t+1}-\y^{t}}_2^2\\&
=\alpha_{t+1}\left(\left<\nabla\psi(\p^{t+1}),\y^{t+1}-\y^{t}\right>+\left<\nabla\psi(\p^{t+1}),\y^{t}-\p^{t+1}\right> - g(\y^{t+1})\right)-\dfrac{1}{2}\norm{\y^{t+1}-\y^{t}}_2^2\\&
\overset{~\eqref{eq:PDDef1},~\eqref{eq:PDDef2}}{=}A_{t}\left<\nabla\psi(\p^{t+1}),\w^{t}-\p^{t+1}\right>-\alpha_{t+1}g(\y^{t+1})+A_{t+1}\left<\nabla\psi(\p^{t+1}),\p^{t+1}-\w^{t+1}\right>\\&-\dfrac{A_{t+1}^2}{2\alpha^{2}_{t+1}}\norm{\p^{t+1}-\w^{t+1}}_2^2\\&
\overset{~\eqref{PDalpQuadEq}}{\leqslant} A_{t}\left<\nabla\psi(\p^{t+1}),\w^{t}-\p^{t+1}\right>-\alpha_{t+1}g(\y^{t+1})\\&+A_{t+1}\left(\left<\nabla\psi(\p^{t+1}),\p^{t+1}-\w^{t+1}\right>-\dfrac{L}{2}\norm{\p^{t+1}-\w^{t+1}}_2^2\right)\\&
\overset{~\eqref{Lip_cond}}{\leqslant} A_{t}\psi(\w^{t})-A_{t}\psi(\p^{t+1})-\alpha_{t+1}g(\y^{t+1})+A_{t+1}\psi(\p^{t+1})-A_{t+1}\psi(\w^{t+1})\\&
=A_{t}\psi(\w^{t})-\alpha_{t+1}g(\y^{t+1})+\alpha_{t+1}\psi(\p^{t+1})-A_{t+1}\psi(\w^{t+1})+A_{t}g(\w^{t})-A_{t}g(\w^{t})\\&
\overset{~\eqref{eq:PDDef2},conv-ty}{\leqslant} A_{t}\psi(\w^{t})+\alpha_{t+1}\psi(\p^{t+1})-A_{t+1}\psi(\w^{t+1})+A_{t}g(\w^{t})-A_{t+1}g(\w^{t+1})\\&
=A_{t}\varphi(\w^{t})+\alpha_{t+1}\psi(\p^{t+1})-A_{t+1}\varphi(\w^{t+1})
\end{split}
\end{equation*}
Finally, using this, we obtain that
\begin{equation*}
\begin{split}
A_{t+1}\varphi(\w^{t+1})-A_{t}\varphi(\w^{t})&\leqslant \alpha_{t+1}\left( \psi(\p^{t+1})+\left<\nabla\psi(\p^{t+1}),\p-\p^{t+1}\right> +g(\p)\right)\\&+\dfrac{1}{2}\norm{\p-\y^t}_2^2
-\dfrac{1}{2}\norm{\p-\y^{t+1}}_2^2.
\end{split}
\end{equation*}
Sum all the inequalities for $t= 0,\; \ldots\; ,N-1$. Then, for any $\p\geqslant \mathtt{0}$
\begin{equation*}
\begin{split}
A_{N}\varphi(\w^{N})-A_{0}\varphi(\w^{0})&\leqslant \sum\limits_{t=0}^{N-1}\alpha_{t+1}\left( \psi(\p^{t+1})+\left<\nabla\psi(\p^{t+1}),\p-\p^{t+1}\right> +g(\p)\right)\\&+\dfrac{1}{2}\norm{\p-\y^0}_2^2
-\dfrac{1}{2}\norm{\p-\y^{N}}_2^2.
\end{split}
\end{equation*}
Whence, since $C_{0}=\alpha_{0}=0$ 
$$
A_{N}\varphi(\w^{N})\leqslant \sum\limits_{t=0}^{N}\alpha_{t}\left( \psi(\p^{t})+\left<\nabla\psi(\p^{t}),\p-\p^{t}\right> +g(\p)\right)+\dfrac{1}{2}\norm{\p-\y^0}_2^2.
$$ 
Taking in the right hand side the minimum in $\p \geqslant \mathtt{0}$ we obtain
\begin{equation}
\label{fin_est}
A_{N}\varphi(\w^{N})\leqslant \min\limits_{p \geqslant \mathtt{0}} \left\{ \sum\limits_{t=0}^{N}\alpha_{t}\left( \psi(\p^{t})+\left<\nabla\psi(\p^{t}),\p-\p^{t}\right> +g(\p)\right)+\dfrac{1}{2}\norm{\p-\y^0}_2^2\right\}.    
\end{equation}
Put the vector $\mathtt{y}^0=(y_1^0, \, \ldots, \,y_n^0)^{\top}$, where $\mathtt{y}^0$~--- vector of initial prices. And its components such as 
$$
0 \leqslant y_k^0 \leqslant p_{max},\; k=1, \, \ldots, \, n,
$$
then we obtain that 
\begin{equation*}
   \norm{\y^{0}}_{2}\leqslant \sqrt{n}p_{max}.
\end{equation*}
Let us introduce a set $B^{+}_{2R}(\mathtt{0})=\left\{ \p \, :\, \p \geqslant 0, \, \norm{\p}_2\leqslant 2R \right\}$, where $R$ determines from~\eqref{4_R}. Since~\eqref{bound_2R} we obtain 
\begin{equation*}
\begin{split}
    A_{N}\varphi(\w^{N})& \leqslant \min\limits_{\p \geqslant \mathtt{0}} \left\{ \sum\limits_{t=0}^{N}\alpha_{t}\left( \psi(\p^{t})+\left<\nabla\psi(\p^{t}),\p-\p^{t}\right> +g(\p)\right)+\dfrac{1}{2}\norm{\p-\y^0}_2^2\right\}\\&
    \leqslant
    \min\limits_{\p \in B_{2R}^{+}(\mathtt{0})}\left\{\sum\limits_{t=0}^{N}\alpha_{t}\left( \psi(\p^{t})+\left<\nabla\psi(\p^{t}),\p-\p^{t}\right> +g(\p)\right)+\dfrac{1}{2}\norm{\p-\y^0}_2^2 \right\}.
\end{split}
\end{equation*}
And since 
$$
\norm{\p-\y^0}_2^2 \leqslant 2\norm{\y^0}_2^2 + 2 \norm{\p}_2^2 \leqslant 8R^2 + \dfrac{2}{9}R^2 = \dfrac{74}{9}R^2
$$
obtain the following
\begin{equation*}
\begin{split}
    A_{N}\varphi(\w^{N})& \leqslant
    \min\limits_{\p \in B_{2R}^{+}(\mathtt{0})}\left\{\sum\limits_{t=0}^{N}\alpha_{t}\left( \psi(\p^{t})+\left<\nabla\psi(\p^{t}),\p-\p^{t}\right> +g(\p)\right) \right\}+\dfrac{37}{9}R^2.
\end{split}
\end{equation*}
Notice, that 
\begin{equation*}
\begin{split}
   &-\min\limits_{\p \in B_{2R}^{+}(\mathtt{0})}\left\{\sum\limits_{t=0}^{N}\alpha_{t}\left( \psi(\p^{t})+\left<\nabla\psi(\p^{t}),\p-\p^{t}\right> +g(\p)\right) \right\} \\&
   =\max\limits_{\p \in B_{2R}^{+}(\mathtt{0})}\left\{\sum\limits_{t=0}^{N}\alpha_{t}\left(\left<\mathtt{x}(\p^{t}),\p^{t}-\p\right> -g(\p)-\left<\mathtt{x}(\p^{t}),\p^{t}\right>+\sum\limits_{k=1}^{n}f_{k}(x_{k}(\p^{t}))\right) \right\}
   \\&= \sum\limits_{t=0}^{N}\alpha_{t}f(\mathtt{x}(\p^{t})) +\max\limits_{\p \in B_{2R}^{+}(\mathtt{0})}\left\{\sum\limits_{t=0}^{N}\alpha_{t} \left( \left<- \mathtt{x}(\p^t),\p\right> - g(\p)\right) \right\} \\&
   \overset{conv-ty}{\geqslant} A_{t}f(\mathtt{x}^N) + A_{t}\max\limits_{\p \in B_{2R}^{+}(\mathtt{0})}\left\{ - \left< \mathtt{x}^N,\p\right>+C\min\limits_{k=1, \, \ldots \, ,n}p_{k} \right\}\\&
   =A_{t}f(\mathtt{x}^N) + A_{t}\max\limits_{\substack{\p \in B_{2R}^{+}(\mathtt{0})\;\\
\\p_1=\ldots=p_n=p}}\left\{  - \left< \mathtt{x}^{N},\p\right>+C\min\limits_{k=1, \, \ldots \, ,n}p_{k} \right\}\\&=
A_{t}f(\mathtt{x}^N) + A_{t}\max\limits_{\substack{\p \in B_{2R}^{+}(\mathtt{0})\;\\
\\p_1=\ldots=p_n=p}}\left\{Cp- \sum\limits_{k=1}^{n}x^{N}_{k} p \right\}\\&
=A_{t}f(\mathtt{x}^N) + A_{t}\dfrac{2R}{\sqrt{n}}\left[C- \sum\limits_{k=1}^{n}x^{N}_k \right]_{+},
\end{split}
\end{equation*}
where $\mathtt{x}^N=\dfrac{1}{A_{N}}\sum\limits_{t=0}^{N}\alpha_t\mathtt{x}(\p^t)$. Thus, we obtain the following estimation 
\begin{equation}
\label{est1}
    \varphi(\w^{N})+f(\mathtt{x}^N)+\dfrac{2R}{\sqrt{n}}\left[C- \sum\limits_{k=1}^{n}x^{N}_k \right]_{+}\leqslant \dfrac{37R^2}{9A_{N}}. 
\end{equation}
Since $\w^*$ is an optimal solution of Problem $(D)$, we have, for any $\mathtt{x}\geqslant 0$
$$
f(\mathtt{x}^*) \leqslant f(\mathtt{x})-\left< \mathtt{x},\w^*\right>+C\min\limits_{k=1, \, \ldots \, ,n}w^*_{k}.
$$
Using the estimation~\eqref{estim_p^*}, we get
$$
f(\mathtt{x}^{N})-f(\mathtt{x}^*)\geqslant -\dfrac{R}{3\sqrt{n}}\left[C- \sum\limits_{k=1}^{n}x^{N}_k \right]_{+}.
$$
Considering the weak duality $f(\mathtt{x}^*)\geqslant -\varphi(\w^*)$, obtain 
\begin{equation*}
\begin{split}
f(\mathtt{x}^N)+\varphi(\w^N)&=f(\mathtt{x}^N)+f(\mathtt{x}^*)-f(\mathtt{x}^*)+\varphi(\w^*)-\varphi(\w^*)+\varphi(\w^N)\\&
\geqslant f(\mathtt{x}^N)-f(\mathtt{x}^*)\geqslant -\dfrac{R}{3\sqrt{n}}\left[C- \sum\limits_{k=1}^{n}x^{N}_k \right]_{+}.
\end{split}   
\end{equation*}
This and~\eqref{est1} give 
\begin{equation}
\label{estimation_constr}
   \left[C- \sum\limits_{k=1}^{n}x^{N}_k \right]_{+}\leqslant \dfrac{37R\sqrt{n}}{15A_{N}}. 
\end{equation}
And also due to the weak duality, we obtain
\begin{equation}
\label{estimation_f}
  f(\mathtt{x}^{N})-f(\mathtt{x}^*)\leqslant \dfrac{37R^2}{9A_{N}}.
\end{equation}

Let us show that the lower bound for the sequence $A_t$, $t \geqslant 0$ is  
\begin{equation}
A_t \geqslant \frac{(t+1)^2}{4L_{\psi}}, \,\, \, \forall t \geqslant 1
\label{eq:AtGrowth}
\end{equation}	
where $L_{\psi}$ is the Lipschitz constant for the gradient of $\psi$. 

For $t=1$, since $\alpha_0 = 0 $ and $A_1 = \alpha_0+\alpha_1= \alpha_1$, we have from \eqref{PDalpQuadEq} 
	\begin{equation*}
		A_1 = \alpha_1 = \frac{1}{L_{\psi}} \geqslant \frac{1}{L_{\psi}}.
	\end{equation*}

	Let us now assume that \eqref{eq:AtGrowth} holds for some $t \geqslant 1$ and prove that it holds for $t+1$.
	From \eqref{PDalpQuadEq}  we have a quadratic equation for $\alpha_{t+1}$
	\begin{equation*}
	L_{\psi}\alpha_{t+1}^2 - \alpha_{t+1} - A_t = 0.
	\end{equation*}
	Since we need to take the largest root, we obtain,
	\begin{align}
	\alpha_{t+1} & = \frac{1 + \sqrt{\uprule 1 + 4L_{\psi}A_{t}}}{2L_{\psi}} = \frac{1}{2L_{\psi}} + \sqrt{\frac{1}{4L_{\psi}^2} + \frac{A_{t}}{L_{\psi}}} \geqslant 
	\frac{1}{2L_{\psi}} + \sqrt{\frac{A_{t}}{L_{\psi}}} \notag \\
	&\geqslant
	\frac{1}{2L_{\psi}} + \frac{1}{\sqrt{L_{\psi}}}\frac{t+1}{2\sqrt{L_{\psi}}} =
	\frac{t+2}{2L_{\psi}}, \notag
	\end{align}
	where we used the induction assumption that \eqref{eq:AtGrowth} holds for $t$.
	Using the obtained inequality, from \eqref{PDalpQuadEq} and \eqref{eq:AtGrowth} for $t$, we get
	\begin{equation*}
	A_{t+1} = A_t + \alpha_{t+1} \geqslant \frac{(t+1)^2}{4L_{\psi}} + \frac{t+2}{2L_{\psi}} \geqslant \frac{(t+2)^2}{4L_{\psi}}.
	\end{equation*}
So, using this estimation and that $R=3p_{max}\sqrt{n}$ and $L_{\psi}=\dfrac{n}{\mu}$ we obtain the statement of the theorem from~\eqref{estimation_constr} and~\eqref{estimation_f}.

\qed

\subsection{The resource allocation problem in vector case}

In this section we consider the vector case of resource allocation problem. Let us consider the same problem as the problem from Section~\ref{Prob_stat}. But now each producer produces $m$ different products, having its own cost function $f_k(\mathtt{x}_k), \, k=1,\, \ldots ,\,n$ representing the total cost of production of the volume $\mathtt{x}_k \in \R^{m}$ which is the number of tons of products produced by the producer $k$ in one year.  So, we have production matrix $X\in \R^{m\times n}$, where each row represents the vector of production for one product by each producer and each column represents the production for one producer of each product. Let us define $\mathtt{x}_{j \bullet }$ the $j$ row of matrix and  $\mathtt{x}_{\bullet k}$ the $k$ column of matrix. The Center buys product from the producers and chooses its strategy in such a way that the total production volume per year by all producers is not less  $c_{m}$ tons of  product $m$. To do so, the Center needs to find $y_{jk}$ - the volume of product $j$ which is purchased from the producer $k$. And in this case we can write down the following resource allocation problem 
\begin{equation*}
\label{Main_problem_2}
(P_1) \quad \quad \sum \limits_{k=1}^n f_k(\mathtt{x}_{\bullet k}) \rightarrow \min \limits_{\substack{\sum \limits_{k=1}^n y_{jk} \geqslant  c_{j}, \;j = 1, \, \ldots, \, m;\; \\ \mathtt{x}_{\bullet k} \geqslant  \y_{\bullet k}, \y_{\bullet k} \geqslant  \mathtt{0},\; k = 1, \, \ldots, \, n,}}
\end{equation*}
where the cost functions~$f_k(\mathtt{x}_{\bullet k})\; k = 1, \, \ldots, \, n$ are increasing for each variable $x_{j,k}$ and~$\mu$-strongly convex. 

Introducing dual variables $\p_{\bullet k}, \, \, k=1, \, \ldots, \, n$ and using the duality, we obtain 
\begin{equation*}
\begin{split}
\min \limits_{\substack{\sum \limits_{k=1}^n y_{jk} \geqslant  c_{j}, \;j = 1, \, \ldots, \, m;\; \\ \mathtt{x}_{\bullet k} \geqslant  \y_{\bullet k}, \y_{\bullet k} \geqslant  \mathtt{0},\; k = 1, \, \ldots, \, n,}} &\sum \limits_{k=1}^n f_k(\mathtt{x}_{\bullet k}) = \min \limits_{\substack{\sum \limits_{k=1}^n y_{jk} \geqslant  c_j, \; \;j = 1, \, \ldots, \, m;\;\\ y_{\bullet k} \geqslant  \mathtt{0},\;k=1, \, \ldots, \, n}} \Bigl\{\sum \limits_{k=1}^n f_k(\mathtt{x}_{\bullet k})+\sum \limits_{k=1}^n \max \limits_{\p_{\bullet k} \geqslant  \mathtt{0}}\p^{\top}_{\bullet k}(\y_{\bullet k}-\mathtt{x}_{\bullet k})\Bigr\}\\
&=-\min\limits_{\p_{\bullet 1}, \, \ldots, \, \p_{\bullet n}\geqslant  0} \Bigl\{\sum \limits_{k=1}^n \max \limits_{\mathtt{x}_{\bullet k} \geqslant  \mathtt{0}}(\p_{\bullet k}^{\top}\mathtt{x}_{\bullet k}-f_k(\mathtt{x}_{\bullet k}))-\min \limits_{\sum \limits_{k=1}^n y_{jk} \geqslant  c_j; \;  \y_{\bullet k} \geqslant  \mathtt{0}}\sum \limits_{k=1}^{n}\p_{\bullet k}^{\top}\y_{\bullet k}\Bigr\}\\
&=-\min\limits_{\p_{\bullet 1}, \, \ldots, \, \p_{\bullet n}\geqslant  \mathtt{0}} \Bigl\{\sum \limits_{k=1}^n \Big\{ \p_k^{\top}\mathtt{x}_{\bullet k}(\p_{\bullet k})-f_k(\mathtt{x}_{\bullet k}(\p_{\bullet k}))\Big\}-\sum \limits_{j=1}^m c_{j}\min \limits_{k=1, \, \ldots, \, n}p_{jk}\Bigr\},
\end{split}
\end{equation*}
where 
\begin{equation}
\mathtt{x}_{\bullet k}(\p_{\bullet k})=\argmax \limits_{\mathtt{x}_{\bullet k} \geqslant  \mathtt{0}} \Big \{ \p_{\bullet k}^{\top} \mathtt{x}_k-f_k(\mathtt{x}_{\bullet k}) \Big \}, \quad  k =1, \, 2, \, \ldots, \, n.
\end{equation}

Then the dual problem (up to a sign) has the following form
\begin{equation*}
(D_1) \quad \quad \varphi(\p_{\bullet 1}, \, \ldots, \, \p_{\bullet n})=\sum \limits_{k=1}^n \Big\{ \p_{\bullet k}^{\top}\mathtt{x}_{\bullet k}(\p_{\bullet k})-f_k(\mathtt{x}_{\bullet k}(\p_{\bullet k}))\Big\}-\sum \limits_{j=1}^m c_{j}\min \limits_{k=1, \, \ldots, \, n}p_{jk} \rightarrow \min \limits_{\p_{\bullet 1}, \, \ldots, \, \p_{\bullet n}\geqslant  \mathtt{0}.}
\end{equation*}
Let us consider this problem as the composite optimization problem. We can rewrite the dual problem as
\begin{equation*}
\varphi(\p_{\bullet 1}, \, \ldots, \, \p_{ \bullet n})=\psi(\p_{\bullet 1}, \, \ldots, \, \p_{ \bullet n})+g(\p_{\bullet 1}, \, \ldots, \, \p_{ \bullet n}),
\end{equation*}
where
\begin{equation}
\label{psi}
\psi(\p_{\bullet 1}, \, \ldots, \, \p_{ \bullet n})=\sum \limits_{k=1}^n \Big\{ \p_{\bullet k}^{\top}\mathtt{x}_{\bullet k}(\p_{\bullet k})-f_k(\mathtt{x}_{\bullet k}(\p_{\bullet k}))\Big\}
\end{equation}
is convex function, which gradient satisfying Lipschitz condition with constant $L_{\psi}=\dfrac{n}{\mu}$. And convex non-smooth composite function 
\begin{equation*}
g(\p_{\bullet 1}, \, \ldots, \, \p_{ \bullet n})=-\sum \limits_{j=1}^m c_{j}\min \limits_{k=1, \, \ldots, \, n}p_{jk}.
\end{equation*}

We define $ \bar{\mathtt{x}}_{\bullet k }=\left(\dfrac{2c_1}{n}, \, \ldots, \, \dfrac{2c_m}{n}\right)^{\top}, \, k=1, \, \ldots, \, n$ and $\bar\y_{\bullet k}=\left(\dfrac{c_1}{n}, \, \ldots, \, \dfrac{c_m}{n}\right)^{\top}, \, k=1, \, \ldots, \, n$, and similarly to Lemma~\ref{Lm:1}, using the Slater's condition, we obtain the following Lemma. 
\begin{Lm}
\label{Lm:4}
Let the $\p_{j \bullet}^{*}, \, j=1, \, \ldots, \, m$ be a solution to the dual problem $(D_1).$ Then it satisfies the inequality 
\begin{equation*}
   \norm{\p_{j \bullet}^{*}}_{2}\leqslant \sqrt{n}\bar p_{max}.
\end{equation*}
where 
\begin{equation}
\label{eq:p_max_vect}
   \bar p_{max}:= \dfrac{n}{\min\limits_{j=1,\ldots,m}c_j}\left(\sum \limits_{k=1}^n f_k\left(\bar{\mathtt{x}}_{ \bullet k}\right)-\sum \limits_{k=1}^n f_k(\mathtt{0})\right). 
\end{equation}
\end{Lm}

\subsubsection{Composite gradient method}
In this subsection to solve the resource allocation problem in vector case we use composite gradient method from the Section~\ref{comp}. Let us define vectors $\bar{\y}=\left(\y_{\bullet 1}^{\top},\, \ldots, \, \y_{ \bullet n}^{\top}\right)\in \R^{mn}$, $\bar{\p}=\left(\p_{\bullet 1}^{\top},\, \ldots, \, \p_{ \bullet n}^{\top}\right)  \in \R^{mn}$ and $\bar{\mathtt{x}}(\bar{\p})=\left(\mathtt{x}_{\bullet 1}(\p_{\bullet 1})^{\top},\, \ldots, \, \mathtt{x}_{ \bullet n}(\p_{ \bullet n})^{\top}\right)  \in \R^{mn}$, then the step of the method can be rewritten as follows
\begin{equation*}
\begin{split}
				\bar{\y}^{t+1}&=\argmin_{\bar{\p}\geqslant 0}\left\{\left<\bar{\mathtt{x}}(\bar{\p}^{t+1}),\bar{\p}-\tilde{\p}^{t+1}\right>+g(\p_{\bullet 1}, \, \ldots, \, \p_{ \bullet n})+\dfrac{ L_{\psi}}{2}\norm{\bar{\p}-\bar{\y}^t}_{2}^{2}\right\}\\& =\argmin_{\bar{\p}\geqslant 0}\left\{\sum_{k=1}^n\left<\mathtt{x}_{\bullet k}(\p_{\bullet k}^{t+1}),\p_{\bullet k}-\p_{\bullet k}^{t+1}\right>-\sum \limits_{j=1}^m c_{j}\min \limits_{k=1, \, \ldots, \, n}p_{jk}+\dfrac{ L_{\psi}}{2}\sum_{k=1}^n\norm{\p_{\bullet k}-\y_{\bullet k}^t}_{2}^{2}\right\} 
				\\& =\argmin_{\bar{\p}\geqslant 0}\left\{\sum \limits_{j=1}^m \left( \left<\mathtt{x}_{j \bullet }(\p_{j \bullet }^{t+1}),\p_{j \bullet }-\p_{j \bullet }^{t+1}\right>- c_{j}\min \limits_{k=1, \, \ldots, \, n}p_{jk}+\dfrac{ L_{\psi}}{2}\norm{\p_{j \bullet }-\y_{j \bullet }^t}_{2}^{2}\right)\right\},
\end{split}
\end{equation*}
where $\p_{j \bullet }=(p_{j1}, \, \ldots, \, p_{jn})^{\top}$, i.e. from the summation of the producers, we moved to the summation of products. Note that we can divide the step into $m$ independent problems for each product, so for product $j$, we have 
\begin{equation*}
				\y_{j \bullet }^{t+1} =\argmin_{\p_{j \bullet }\geqslant 0}\left\{ \left<\mathtt{x}_{j \bullet }(\p_{j \bullet }^{t+1}),\p_{j \bullet }-\p_{j \bullet }^{t+1}\right>- c_{j}\min \limits_{k=1, \, \ldots, \, n}p_{jk}+\dfrac{\bar L_{\psi}}{2}\norm{\p_{j \bullet }-\y_{j \bullet }^t}_{2}^{2}\right\},
\end{equation*}
where~$\p_{j \bullet }$~-- is price vector of all producers for product $j$ and $\mathtt{x}_{j \bullet }(\p_{j \bullet })$~-- is vector of optimal plans for the production of product $j$ for all producers. Define $\tilde{\p}_{j \bullet }^{t+1}=\p^t_{j \bullet }-\dfrac{1}{L_{\psi}}\mathtt{x}_{j \bullet }(\p_{j \bullet }^{t+1})$ 
and then, using Lemma~\ref{lemma2}, we obtain the following solution:
\begin{itemize}
    \item If $\sum\limits_{k=1}^{n}\left(-\tilde{p}^{t+1}_{jk} \right)_{+} \geqslant \dfrac{c_j}{L_{\psi}} $ then $p^{t+1}_{j.center}=0$. 
    \item Else $p^{t+1}_{j.center} > 0$ is determined from 
    $$
    \sum\limits_{k=1}^{n}\left(p^{t+1}_{j.center}- \tilde{p}^{t+1}_{jk} \right)_{+}=\dfrac{c_j}{L_{\psi}} 
    $$
    \end{itemize}
    and the solution is determined as 
    $$
    p^{t+1}_{jk}=\max\left(p^{t+1}_{j.center},\,  \tilde{p}^{t+1}_{jk}\right), \, \, k=1,\, \ldots\, ,n.
    $$

\begin{Th}
\label{th6}
Let Algorithm \eqref{Alg5} be run for $N$ steps with starting points $\p_{j\bullet}^0, \, \, j=1, \, \ldots, \, m$ satisfying $0 \leqslant p_{jk}^0 \leqslant \bar{p}_{max},\; k=1, \, \ldots, \, n$, where $\bar{p}_{max}$ is given in \eqref{eq:p_max_vect}. Then
\begin{equation*}
\begin{split}
\sum_{k=1}^{n}f_k(\mathtt{x}_{\bullet k}^N)-\sum_{k=1}^{n}f(\mathtt{x}_{\bullet k}^{*})&\leqslant   \sum_{k=1}^{n}f_k(\mathtt{x}_{\bullet k}^N)+\varphi(\bar {\mathtt{p}}^{*})\leqslant   \varphi(\bar {\mathtt{p}}^N)+\sum_{k=1}^{n}f_k(\mathtt{x}_{\bullet k}^N) \leqslant \dfrac{82 \bar p_{max}^2 n^2 m}{N\mu} ,\\ 
 \sum_{j=1}^{m} \left[c_j- \sum\limits_{k=1}^{n}x_{jk}^N \right]_{+} & \leqslant \dfrac{82 \bar p_{max} n^2 m}{3N\mu}, 
\end{split}
\end{equation*}
where $\p^{N}_{\bullet k}=\dfrac{1}{N}\sum\limits_{t=1}^{N}\p_{\bullet k}^{t}$, $\mathtt{x}_{\bullet k}^{N}=\dfrac{1}{N}\sum\limits_{t=0}^{N-1}\mathtt{x}_{\bullet k}(\p_{\bullet k}^{t})$ and $\bar\p^{N}=(\p_{\bullet 1}^{N}, \, \ldots, \, \p_{\bullet n}^{N})$
\end{Th}
\textbf{Proof of Theorem~
\ref{th6}.}
Similarly to the proof of the theorem~\ref{th2} we can obtain the equation~\eqref{fin_est1}, but in this case we can write dawn this equation as 
\begin{equation*}
\varphi(\bar\p^{N})\leqslant \dfrac{1}{N}\min\limits_{\bar \p \geqslant \mathtt{0}} \left\{ \sum\limits_{t=0}^{N}\sum_{j=1}^m\left<\mathtt{x}_{j \bullet}(\p_{j \bullet}^{t+1}),\p_{j \bullet}\right>-\sum_{k=1}^n f_k(\mathtt{x}_{ \bullet k}) +g(\bar \p)+\dfrac{\bar L_{\psi}}{2}\sum_{j=1}^m\norm{\p_{j \bullet}-\p_{j \bullet}^0}_2^2\right\},  \end{equation*}
where $\bar\p^{N}=(\p_1^{N}, \, \ldots, \, \p_n^{N})$. Let us introduce a set $B^{+}_{3\bar R}(\mathtt{0})=\left\{ \p \, :\, \p \geqslant 0, \, \norm{\p}_2\leqslant 2\bar R \right\}$, where $\bar R=3\sqrt{n}\bar p_{max}$ determines similarly to~\eqref{4_R}. Since we also have equation similarly to~\eqref{bound_2R} we obtain that $\p^t_{j \bullet} \in B^{+}_{3\bar R}(\mathtt{0}), \, \forall t.$ So, we obtain 
\begin{equation*}
\begin{split}
\varphi&(\bar\p^{N}) \leqslant \dfrac{1}{N} \min\limits_{\bar \p \geqslant \mathtt{0}} \left\{ \sum\limits_{t=0}^{N}\sum_{j=1}^m\left<\mathtt{x}_{j \bullet}(\p_{j \bullet}^{t+1}),\p_{j \bullet}\right>-\sum_{k=1}^n f_k(\mathtt{x}_{ \bullet k}) +g(\bar \p)+\dfrac{\bar L_{\psi}}{2}\sum_{j=1}^m\norm{\p_{j \bullet}-\p_{j \bullet}^0}_2^2\right\}\\&
=\dfrac{1}{N}\min\limits_{\bar \p_{j \bullet} \in B^{+}_{3\bar R}(\mathtt{0})} \left\{ \sum\limits_{t=0}^{N}\sum_{j=1}^m\left<\mathtt{x}_{j \bullet}(\p_{j \bullet}^{t+1}),\p_{j \bullet}\right>-\sum_{k=1}^n f_k(\mathtt{x}_{ \bullet k}) +g(\bar \p)+\dfrac{\bar L_{\psi}}{2}\sum_{j=1}^m\norm{\p_{j \bullet}-\p_{j \bullet}^0}_2^2\right\}
\end{split}
\end{equation*}
And since 
$$
\norm{\p_{j \bullet}-\p_{j \bullet}^0}_2^2 \leqslant 2\norm{\p_{j \bullet}^0}_2^2 + 2 \norm{\p_{j \bullet}}_2^2 \leqslant 18 \bar R^2 + \dfrac{2}{9}\bar R^2 = \dfrac{164}{9}\bar R^2
$$
obtain the following
\begin{equation*}
\varphi(\bar\p^{N}) \leqslant \dfrac{1}{N} \min\limits_{\bar \p_{j \bullet} \in B^{+}_{3\bar R}(\mathtt{0})} \left\{ \sum\limits_{t=0}^{N}\sum_{j=1}^m\left<\mathtt{x}_{j \bullet}(\p_{j \bullet}^{t+1}),\p_{j \bullet}\right>-\sum_{k=1}^n f_k(\mathtt{x}_{\bullet k}) +g(\bar \p)\right\}+\dfrac{82m\bar L_{\psi}}{9N}\bar R^2.
\end{equation*}
Notice, that 
\begin{equation*}
\begin{split}
    &-\dfrac{1}{N}\min\limits_{\bar \p_{j \bullet} \in B^{+}_{3\bar R}(\mathtt{0})} \left\{ \sum\limits_{t=0}^{N}\left(\sum_{j=1}^m\left<\mathtt{x}_{j \bullet}(\p_{j \bullet}^{t+1}),\p_{j \bullet}\right>-\sum_{k=1}^n f_k(\mathtt{x}_{\bullet k}) +g(\bar \p)\right)\right\}\\&=\dfrac{1}{N}\sum_{t=0}^{N}\sum_{k=1}^{n}f_k(\mathtt{x}_{\bullet k})+\dfrac{1}{N}
    \max\limits_{\bar \p_{j \bullet} \in B^{+}_{3\bar R}(\mathtt{0})} \left\{ \sum\limits_{t=0}^{N}-\sum_{j=1}^m\left<\mathtt{x}_{j \bullet}(\p_{j \bullet}^{t+1}),\p_{j \bullet}\right> +\sum \limits_{j=1}^m c_{j}\min \limits_{k=1, \, \ldots, \, n}p_{jk} \right\}\\&
    \overset{conv-ty}{\geqslant} \sum_{k=1}^n f_k(\mathtt{x}^N_{\bullet k})+\max\limits_{\substack{\p_{j \bullet} \in B_{3\bar R}^{+}(\mathtt{0})\;\\
\\p_{j1}=\ldots=p_{jn}=p_j}} \left\{ \sum_{j=1}^m\left(-\left<\mathtt{x}_{j \bullet}^N,\p_{j \bullet}\right> +c_{j}\min \limits_{k=1, \, \ldots, \, n}p_{jk}\right) \right\}\\&=
 \sum_{k=1}^n f_k(\mathtt{x}^N_{\bullet k})+\max\limits_{\substack{\p_{j \bullet} \in B_{3\bar R}^{+}(\mathtt{0})\;\\
\\p_{j1}=\ldots=p_{jn}=p_j}} \left\{ \sum_{j=1}^m\left(-\sum_{k=1}^n x^N_{jk}p_j+c_{j}p_{j}\right) \right\}\\&= \sum_{k=1}^n f_k(\mathtt{x}^N_{\bullet k})+\dfrac{3\bar R}{\sqrt{n}}\sum_{j=1}^m\left[c_{j}-\sum_{k=1}^n x^N_{jk}\right]_{+},
\end{split}
\end{equation*}
where $\mathtt{x}_{ \bullet k}^N=\dfrac{1}{N}\sum\limits_{t=0}^{N}\mathtt{x}_{ \bullet k}(\p_{\bullet k}^t)$. Thus, we obtain the following estimation 
\begin{equation*}
    \varphi(\bar \w^{N})+\sum_{k=1}^n f_k(\mathtt{x}^N_{\bullet k})+\dfrac{3\bar R}{\sqrt{n}}\sum_{j=1}^m\left[c_{j}-\sum_{k=1}^n x^N_{jk}\right]_{+}\leqslant \dfrac{82 \bar L_{\psi}\bar R^2m}{9N}. 
\end{equation*}
And similarly to~\eqref{duality_gap1} we obtain 
$$
\sum_{j=1}^m\left[c_{j}-\sum_{k=1}^n x^N_{jk}\right]_{+} \leqslant \dfrac{82 \bar L_{\psi}\bar Rm\sqrt{n}}{9N}.
$$
And due to the weak duality, we obtain 
$$
\sum_{k=1}^n f_k(\mathtt{x}^N_{\bullet k})-\sum_{k=1}^n f_k(\mathtt{x}^*_{\bullet k})\leqslant \dfrac{82 \bar L_{\psi}\bar R^2m}{9N}.
$$
And since $\bar{R}=3\sqrt{n}\bar{p}_{max}$ and $L_{\psi}=\dfrac{n}{\mu}$ we obtain the statement of the theorem.

\qed

\begin{equation}
\label{Alg5}
\begin{array}{|c|}
\hline \\
\mbox{\bf Composite gradient method for the resource allocation 
(vector case)}\\
\\
\hline \\
\quad \mbox{
\begin{minipage}{13cm}
\textbf{Input:} $N> 0$ -- number of steps, $L_{\psi}$ -- Lipschitz constant of gradient $\psi$, $\p_{j\bullet}^0, \, \, j=1, \, \ldots, \, m$ -- starting point.
\begin{enumerate}
	\item  Knowing the prices $p_{\bullet k}^t, \, k=1, \, \ldots, \, n$ for the current year $t$, producers calculate  the optimal plan for the production according these prices as 
	\begin{equation*}
\mathtt{x}_{\bullet k}(\p^{t+1}_{\bullet k})=\argmax \limits_{\mathtt{x}_k \geqslant  0} \Big \{ \sum_{j=1}^m p^{t+1}_{jk} x_{jk}-f_k(\mathtt{x}_k) \Big \}, \quad  k =1, \, 2, \, \ldots, \, n.
\end{equation*}
\item Each factory predicts the price for the next year $t+1$ for product $j=1,\, \ldots\, ,m$ as 
$$
\tilde{p}^{t+1}_{jk}=p^t_{jk}-\dfrac{1}{L_{\psi}}x_{jk}(p_{jk}^t), \quad  k =1, \, 2, \, \ldots, \, n
$$
and send this information to the Center.

	\item The Center determines the price $p^{t+1}_{j.center}$ at which it will purchase product $j$ for the next year $t+1$ as 
\begin{itemize}
    \item If $\sum\limits_{k=1}^{n}\left(-\tilde{p}^{t+1}_{jk} \right)_{+} \geqslant \dfrac{C}{ L_{\psi}}$ then $p^{t+1}_{j.center}=0$;
    \item Else $p^{t+1}_{j.center} > 0$ is determined from 
    $$
    \sum\limits_{k=1}^{n}\left(p^{t+1}_{j.center}-\tilde{p}^{t+1}_{jk} \right)_{+}= \dfrac{C}{ L_{\psi}},
    $$
\end{itemize}
and sends this price to all producers.
	
	\item  Each producer adjusts the price for product $ j=1,\, \ldots\, ,m$ for the next year as follows
$$
    p^{t+1}_{jk}=\max\left(p^{t+1}_{j.center},\, \tilde{p}^{t+1}_{jk}\right), \, \, k=1,\, \ldots\, ,n.
    $$

\end{enumerate}
\end{minipage}
}
\quad\\
\\
\hline
\end{array}
\end{equation}

\subsubsection{Accelerated composite gradient method}
In this subsection, to solve the resource allocation problem in vector case, we use accelerated composite gradient method from the Section~\ref{AGM}. Similarly to the previous subsection the step of the method we can be rewritten as
\begin{equation*}
\begin{split}
				\bar{\y}^{t+1}&=\argmin_{\bar{\p}\geqslant 0}\left\{\alpha_{t+1}\left(\left<\bar{\mathtt{x}}(\bar{\p}^{t+1}),\bar{\p}-\tilde{\p}^{t+1}\right>+g(\p_{\bullet 1}, \, \ldots, \, \p_{ \bullet n})\right)+\dfrac{1}{2}\norm{\bar{\p}-\bar{\y}^t}_{2}^{2}\right\}\\& =\argmin_{\bar{\p}\geqslant 0}\left\{\alpha_{t+1}\left(\sum_{k=1}^n\left<\mathtt{x}_{\bullet k}(\p_{\bullet k}^{t+1}),\p_{\bullet k}-\p_{\bullet k}^{t+1}\right>-\sum \limits_{j=1}^m c_{j}\min \limits_{k=1, \, \ldots, \, n}p_{jk}\right)+\dfrac{1}{2}\sum_{k=1}^n\norm{\p_{\bullet k}-\y_{\bullet k}^t}_{2}^{2}\right\} 
				\\& =\argmin_{\bar{\p}\geqslant 0}\left\{\sum \limits_{j=1}^m \left( \alpha_{t+1}\left(\left<\mathtt{x}_{j \bullet }(\p_{j \bullet }^{t+1}),\p_{j \bullet }-\p_{j \bullet }^{t+1}\right>- c_{j}\min \limits_{k=1, \, \ldots, \, n}p_{jk}\right)+\dfrac{1}{2}\norm{\p_{j \bullet }-\y_{j \bullet }^t}_{2}^{2}\right)\right\},
\end{split}
\end{equation*}
where $\p_{j \bullet }=(p_{j1}, \, \ldots, \, p_{jn})^{\top}$, i.e. from the summation of the producers, we moved to the summation of products. Note that we can divide the step into $m$ independent problems for each product, so for product $j$ we have 
\begin{equation*}
				\y_{j \bullet }^{t+1} =\argmin_{\p_{j \bullet }\geqslant 0}\left\{ \alpha_{t+1}\left(\left<\mathtt{x}_{j \bullet }(\p_{j \bullet }^{t+1}),\p_{j \bullet }-\p_{j \bullet }^{t+1}\right>- c_{j}\min \limits_{k=1, \, \ldots, \, n}p_{jk}\right)+\dfrac{1}{2}\norm{\p_{j \bullet }-\y_{j \bullet }^t}_{2}^{2}\right\},
\end{equation*}
where~$\p_{j \bullet }$~-- is price vector of all producers for product $j$ and $\mathtt{x}_{j \bullet }(\p_{j \bullet })$~-- is vector of optimal plans for the production of product $j$ for all producers.
Define $\tilde{\y}_{j \bullet }^{t+1}=\y^t_{j \bullet }-\alpha_{t+1}\mathtt{x}_{j \bullet }(\p_{j \bullet }^{t+1})$ 
and then, using Lemma~\ref{lemma2}, we obtain the following solution:
\begin{itemize}
    \item If $\sum\limits_{k=1}^{n}\left(-\tilde{y}^{t+1}_{jk} \right)_{+}\geqslant c_{j}\alpha_{t+1}$ then $y^{t+1}_{j.center}=0$. 
    \item Else $y^{t+1}_{j.center} > 0$ is determines from 
    $$
    \sum\limits_{k=1}^{n}\left(y^{t+1}_{j.center}- \tilde{y}^{t+1}_{jk} \right)_{+}= c_j\alpha_{t+1}
    $$
    \end{itemize}
    and the solution is determined as 
    $$
    y^{t+1}_{jk}=\max\left(y^{t+1}_{j.center},\,  \tilde{y}^{t+1}_{jk}\right), \, \, k=1,\, \ldots\, ,n.
    $$

\begin{equation}
\label{Alg6}
\begin{array}{|c|}
\hline \\
\mbox{\bf Accelerated composite gradient descent for resource allocation
(vector case)}\\
\\
\hline \\
\quad \mbox{
\begin{minipage}{13cm}
\textbf{Input:} $N> 0$ -- number of steps, $L_{\psi}$ -- Lipschitz constant of gradient $\psi$, $\p_{j\bullet}^0=\y_{j\bullet}^0=\w_{j\bullet}^0, \, \, j=1, \, \ldots, \, m$ -- starting points.
\begin{enumerate}
	\item  In the current year producers find $\alpha_{t+1}$ as the largest root of the equation
				\begin{equation*}
				A_{t+1}:=A_t+\alpha_{t+1} =  L_{\psi}\alpha_{t+1}^2.
				\end{equation*}

		\item All producers calculate the average price $\p_{j \bullet }^{t+1}$ for product $j=1,\, \ldots\, ,m$ as 
				\begin{equation*}
				p_{jk}^{t+1} = \frac{\alpha_{t+1}y^t_{jk} + A_t w_{jk}^t}{A_{t+1}}, \, \, k=1,\, \ldots\, ,n 
				\end{equation*}
	 and calculate  the optimal plan for the production as 
	\begin{equation*}
\mathtt{x}_{\bullet k}(\p^{t+1}_{\bullet k})=\argmax \limits_{\mathtt{x}_k \geqslant  0} \Big \{ \sum_{j=1}^m p^{t+1}_{jk} x_{jk}-f_k(\mathtt{x}_k) \Big \}, \quad  k =1, \, 2, \, \ldots, \, n.
\end{equation*}
\item Each producer predict the price for the next year $t+1$ for product $j=1,\, \ldots\, ,m$ as 
$$
\tilde{y}_{jk}^{t+1}=y^t_{jk}-\alpha_{t+1}x_{jk}(\p_k^{t+1}), \, \, k=1,\, \ldots\, ,n
$$
and send this information to the Center.
	\item The Center determines the prediction prices $y^{t+1}_{j.center}$  for each product $j=1,\, \ldots\, ,m$ for the next year  as 
\begin{itemize}
    \item If $\sum\limits_{k=1}^{n}\left(-\tilde{y}^{t+1}_{jk} \right)_{+} \geqslant c_j\alpha_{t+1}$ then $y^{t+1}_{j.center}=0$ 
    \item Else $y^{t+1}_{j.center} > 0$ and determines from 
    $$
    \sum\limits_{k=1}^{n}\left(y^{t+1}_{j.center}- \tilde{y}^{t+1}_{jk} \right)_{+}= c_j\alpha_{t+1}.
    $$
    \end{itemize}
    \item  Each producer adjusts the prediction price for the next year as follows
 $$
    y^{t+1}_{jk}=\max\left(y^{t+1}_{j.center},\,  \tilde{y}^{t+1}_{jk}\right), \, \, k=1,\, \ldots\, ,n.
    $$
and calculate the historical price for product $j=1,\, \ldots\, ,m$ for the next year 
				\begin{equation*}
				w_{jk}^{t+1} = \frac{\alpha_{t+1}y_{jk}^{t+1} + A_t w_{jk}^t}{A_{t+1}}, \, \, k=1,\, \ldots\, ,n.
				\end{equation*}
\end{enumerate}
\end{minipage}
}
\quad\\
\\
\hline
\end{array}
\end{equation}

\begin{Th}
\label{th5}
Let Algorithm \eqref{Alg6} be run for $N$ steps with starting points $\p_{j\bullet}^0=\y_{j\bullet}^0=\w_{j\bullet}^0, \, \, j=1, \, \ldots, \, m$ satisfying $0 \leqslant p_{jk}^0 \leqslant \bar{p}_{max},\; k=1, \, \ldots, \, n$, where $\bar{p}_{max}$ is given in \eqref{eq:p_max_vect}. Then
\begin{equation*}
\begin{split}
\sum_{k=1}^{n}f_k(\mathtt{x}_{\bullet k}^N)-\sum_{k=1}^{n}f(\mathtt{x}_{\bullet k}^{*})&\leqslant   \sum_{k=1}^{n}f_k(\mathtt{x}_{\bullet k}^N)+\varphi(\bar {\mathtt{w}}^{*})\leqslant   \varphi(\bar {\mathtt{w}}^N)+\sum_{k=1}^{n}f_k(\mathtt{x}_{\bullet k}^N) \leqslant \dfrac{148\bar p_{max}^2n^2 m}{(N+1)^2\mu} ,\\ 
 \sum_{j=1}^{m} \left[c_j- \sum\limits_{k=1}^{n}x_{jk}^N \right]_{+} & \leqslant \dfrac{148\bar p_{max}mn^2}{5(N+1)^2\mu},
\end{split}
\end{equation*}
where $\mathtt{x}_{\bullet k}^N=\dfrac{1}{A_{N}}\sum\limits_{t=0}^{N}\alpha_t\mathtt{x}_{\bullet k}(\p_{\bullet k}^t)$ and $\bar\w^{N}=(\w_{\bullet 1}^{N}, \, \ldots, \, \w_{\bullet n}^{N})$.
\end{Th}

\textbf{Proof of Theorem~
\ref{th5}.}
Similarly to the proof of the theorem~\ref{th3} we can obtain the equation~\eqref{fin_est}, but in this case we can write dawn this equation as 
\begin{equation*}
A_{N}\varphi(\bar\w^{N})\leqslant \min\limits_{\bar \p \geqslant \mathtt{0}} \left\{ \sum\limits_{t=0}^{N}\alpha_{t}\left(\sum_{j=1}^m\left<\mathtt{x}_{j \bullet}(\p_{j \bullet}^{t+1}),\p_{j \bullet}\right>-\sum_{k=1}^n f_k(\mathtt{x}_{ \bullet k}) +g(\p)\right)+\dfrac{1}{2}\sum_{j=1}^m\norm{\p_{j \bullet}-\y_{j \bullet}^0}_2^2\right\},  \end{equation*}
where $\bar\w^{N}=(\w_1^{N}, \, \ldots, \, \w_n^{N})$.  Since we also have equations similarly to~\eqref{bound_2R} we obtain that $\p^t_{j \bullet},\, \w^t_{j \bullet}, \, \y^t_{j \bullet} \in B^{+}_{2\bar R}(\mathtt{0}), \, \forall t, $ where $B^{+}_{2\bar R}(\mathtt{0})=\left\{ \p \, :\, \p \geqslant 0, \, \norm{\p}_2\leqslant 2\bar R \right\}$, $\bar R=3\sqrt{n}\bar p_{max}$. So, we obtain 
\begin{equation*}
\begin{split}
A_{N}\varphi&(\bar\w^{N}) \leqslant \min\limits_{\bar \p \geqslant \mathtt{0}} \left\{ \sum\limits_{t=0}^{N}\alpha_{t}\left(\sum_{j=1}^m\left<\mathtt{x}_{j \bullet}(\p_{j \bullet}^{t+1}),\p_{j \bullet}\right>-\sum_{k=1}^n f_k(\mathtt{x}_{ \bullet k}) +g(\p)\right)+\dfrac{1}{2}\sum_{j=1}^m\norm{\p_{j \bullet}-\y_{j \bullet}^0}_2^2\right\}\\&
=\min\limits_{\bar \p_{j \bullet} \in B^{+}_{2\bar R}(\mathtt{0})} \left\{ \sum\limits_{t=0}^{N}\alpha_{t}\left(\sum_{j=1}^m\left<\mathtt{x}_{j \bullet}(\p_{j \bullet}^{t+1}),\p_{j \bullet}\right>-\sum_{k=1}^n f_k(\mathtt{x}_{ \bullet k}) +g(\p)\right)+\dfrac{1}{2}\sum_{j=1}^m\norm{\p_{j \bullet}-\y_{j \bullet}^0}_2^2\right\}
\end{split}
\end{equation*}
And since 
$$
\norm{\p_{j \bullet}-\y_{j \bullet}^0}_2^2 \leqslant 2\norm{\y_{j \bullet}^0}_2^2 + 2 \norm{\p_{j \bullet}}_2^2 \leqslant 8\bar R^2 + \dfrac{2}{9}\bar R^2 = \dfrac{74}{9}\bar R^2
$$
obtain the following
\begin{equation*}
A_{N}\varphi(\bar\w^{N}) \leqslant \min\limits_{\bar \p_{j \bullet} \in B^{+}_{2\bar R}(\mathtt{0})} \left\{ \sum\limits_{t=0}^{N}\alpha_{t}\left(\sum_{j=1}^m\left<\mathtt{x}_{j \bullet}(\p_{j \bullet}^{t+1}),\p_{j \bullet}\right>-\sum_{k=1}^n f_k(\mathtt{x}_{\bullet k}) +g(\p)\right)\right\}+\dfrac{37m}{9}\bar R^2.
\end{equation*}
Notice, that 
\begin{eqnarray*}
    &-&\min\limits_{\bar \p_{j \bullet} \in B^{+}_{2\bar R}(\mathtt{0})} \left\{ \sum\limits_{t=0}^{N}\alpha_{t}\left(\sum_{j=1}^m\left<\mathtt{x}_{j \bullet}(\p_{j \bullet}^{t+1}),\p_{j \bullet}\right>-\sum_{k=1}^n f_k(\mathtt{x}_{\bullet k}) +g(\p)\right)\right\}\\&=&\sum_{t=0}^{N}\sum_{k=1}^{n}f_k(\mathtt{x}_{\bullet k})+
    \max\limits_{\bar \p_{j \bullet} \in B^{+}_{2\bar R}(\mathtt{0})} \left\{ \sum\limits_{t=0}^{N}\alpha_{t}\left(-\sum_{j=1}^m\left<\mathtt{x}_{j \bullet}(\p_{j \bullet}^{t+1}),\p_{j \bullet}\right> +\sum \limits_{j=1}^m c_{j}\min \limits_{k=1, \, \ldots, \, n}p_{jk} \right)\right\}\\&
    \overset{conv-ty}{\geqslant}& A_t\sum_{k=1}^n f_k(\mathtt{x}^N_{\bullet k})+A_t\max\limits_{\substack{\p_{j \bullet} \in B_{2\bar R}^{+}(\mathtt{0})\;\\
\\p_{j1}=\ldots=p_{jn}=p_j}} \left\{ \sum_{j=1}^m\left(-\left<\mathtt{x}_{j \bullet}^N,\p_{j \bullet}\right> +c_{j}\min \limits_{k=1, \, \ldots, \, n}p_{jk}\right) \right\}\\&=&
 A_t\sum_{k=1}^n f_k(\mathtt{x}^N_{\bullet k})+A_t\max\limits_{\substack{\p_{j \bullet} \in B_{2\bar R}^{+}(\mathtt{0})\;\\
\\p_{j1}=\ldots=p_{jn}=p_j}} \left\{ \sum_{j=1}^m\left(-\sum_{k=1}^n x^N_{jk}p_j+c_{j}p_{j}\right) \right\}\\&=& A_t\sum_{k=1}^n f_k(\mathtt{x}^N_{\bullet k})+A_t\dfrac{2\bar R}{\sqrt{n}}\sum_{j=1}^m\left[c_{j}-\sum_{k=1}^n x^N_{jk}\right]_{+},
\end{eqnarray*}
where $\mathtt{x}_{ \bullet k}^N=\dfrac{1}{A_{N}}\sum\limits_{t=0}^{N}\alpha_t\mathtt{x}_{ \bullet k}(\p_{\bullet k}^t)$. Thus, we obtain the following estimation 
\begin{equation*}
    \varphi(\bar \w^{N})+\sum_{k=1}^n f_k(\mathtt{x}^N_{\bullet k})+\dfrac{2\bar R}{\sqrt{n}}\sum_{j=1}^m\left[c_{j}-\sum_{k=1}^n x^N_{jk}\right]_{+}\leqslant \dfrac{37\bar R^2m}{9A_{N}}. 
\end{equation*}
And similarly to~\eqref{estimation_constr} we obtain 
$$
\sum_{j=1}^m\left[c_{j}-\sum_{k=1}^n x^N_{jk}\right]_{+} \leqslant \dfrac{37\bar Rm}{15A_{N}}.
$$
And due to the weak duality, we obtain 
$$
\sum_{k=1}^n f_k(\mathtt{x}^N_{\bullet k})-\sum_{k=1}^n f_k(\mathtt{x}^*_{\bullet k})\leqslant \dfrac{37\bar R^2m}{9A_{N}}.
$$
And using the estimation~\eqref{eq:AtGrowth} for the $A_t, \, t \geqslant 0$ we obtain the statement on the Theorem.

\qed

\end{document}